\documentclass[a4paper]{article}

\usepackage{amsfonts}
\usepackage{amsmath}
\usepackage{graphicx}
\usepackage{epstopdf}
\usepackage{algorithmic}
\usepackage{tikz}
\usepackage{pgfplots}
\usepackage{pgfplotstable}
\usetikzlibrary{positioning}
\makeatletter
\definecolor{gnuplot@orange}{RGB}{229,158,0}
\definecolor{gnuplot@purple}{RGB}{148,0,212}
\definecolor{gnuplot@lightblue}{RGB}{87,181,232}
\definecolor{gnuplot@green}{RGB}{0,158,115}
\definecolor{gnuplot@darkblue}{RGB}{0,115,179}
\definecolor{gnuplot@yellow}{RGB}{240,227,66}
\makeatother

\ifpdf
  \DeclareGraphicsExtensions{.eps,.pdf,.png,.jpg}
\else
  \DeclareGraphicsExtensions{.eps}
\fi

\newcommand{\TheTitle}{A performance comparison of continuous and
  discontinuous Galerkin methods with fast multigrid solvers}

\title{{\TheTitle}\thanks{This
    work was partially supported by the German Research Foundation (DFG) under
    the project ``High-order discontinuous Galerkin for the exa-scale''
    (ExaDG) within the priority program ``Software for Exascale Computing''
    (SPPEXA), grant agreement no. KR4661/2-1 and WA1521/18-1. The authors gratefully acknowledge the Gauss Centre for Supercomputing e.V.~(\texttt{www.gauss-centre.eu}) for funding this project
by providing computing time on the GCS Supercomputer SuperMUC at Leibniz Supercomputing Centre (LRZ, \texttt{www.lrz.de})
through project id pr83te.}}

\author{
  Martin Kronbichler\thanks{Institute for Computational Mechanics, Technical University of Munich, Boltzmannstr.~15, 85748 Garching b.~M\"unchen, Germany
    (\texttt{\{kronbichler,wall\}@lnm.mw.tum.de}).}
  \and
  Wolfgang A. Wall\footnotemark[2]
}

\usepackage{amsopn}

\begin{document}

\maketitle

\begin{abstract}
  This study presents a fair performance comparison of the continuous finite element method, the symmetric interior penalty discontinuous Galerkin method, and the hybridized discontinuous Galerkin method. Modern implementations of high-order methods with state-of-the-art multigrid solvers for the Poisson equation are considered, including fast matrix-free implementations with sum factorization on quadrilateral and hexahedral elements. For the hybridized discontinuous Galerkin method, a multigrid approach that combines a grid transfer from the trace space to the space of linear finite elements with algebraic multigrid on further levels is developed. Despite similar solver complexity of the matrix-based HDG solver and matrix-free geometric multigrid schemes with continuous and discontinuous Galerkin finite elements, the latter offer up to order of magnitude faster time to solution, even after including the superconvergence effects. This difference is because of vastly better performance of matrix-free operator evaluation as compared to sparse matrix-vector products. A roofline performance model confirms the advantage of the matrix-free implementation.
\end{abstract}

\section{Introduction}
\label{intro}

The relative efficiency of various realizations of the discontinuous Galerkin (DG) method as compared to
continuous finite elements (continuous Galerkin, CG) has been the subject of a
number of recent studies \cite{gmfh13,harp13,ksc12,Yakovlev16}. The
hybridizable discontinuous Galerkin (HDG) method \cite{cgl09,npc09} has
attracted particular interest because it promises a more efficient solution of
linear systems than other discontinuous Galerkin methods in terms of the
number of degrees of freedom and nonzero entries in the system matrix. As
opposed to continuous finite elements or symmetric interior penalty
discontinuous Galerkin \cite{abcm02} methods that rely on the primal
formulation of the differential equation, the HDG method poses the problem in
mixed form using an additional variable for the flux. In order to avoid
solving a linear system involving both the primal and flux unknowns, an
initially counter-intuitive step is taken by introducing yet another variable,
the so-called trace variable defined on the mesh skeleton. The numerical
fluxes in the mixed system are solely expressed in terms of the trace variable
and the local unknowns, avoiding direct coupling between neighboring
elements. As a consequence,
all element unknowns of the primal and flux variables can be eliminated prior to
solving the global linear system by an element-by-element Schur complement,
resulting in a global system in terms of the trace variable only. This Schur complement
approach is conceptually the same as the technique of static condensation in
continuous finite elements that eliminates the unknowns with coupling inside a
single element as a means for improving the efficiency of the solution stage,
see e.g.~\cite{cgl09,gmfh13,harp13,loehner13,Yakovlev16} and references
therein.

Previous efficiency comparisons have found that HDG is a highly competitive
option for two-dimensional problems and direct solvers
\cite{harp13,ksc12}. This work widens the perspective by considering
large-scale problems which demand for iterative solvers and optimal complexity
preconditioners. In that setting it is not enough to characterize the
sparsity structure in the linear system of equations or the nonzero entries in
the matrix as a proxy for the cost of one operator evaluation. Instead, the
interesting factors for competitive solver times are the preconditioner
efficiency and complexity, i.e., the iteration counts, the number of matrix-vector products
per iteration, and timings for one matrix-vector product. This work considers
multigrid methods which are among the most competitive solvers for
elliptic operators on general meshes \cite{Gholami16}. Geometric multigrid (GMG) methods
combine simple iterative schemes on a hierarchy of coarser meshes.  Different
error frequencies are attacked on different mesh levels, such that simple
iterative schemes that smooth the respective high
frequencies on each level can be used. For certain applications, GMG is too restrictive because
a mesh coarsening must be explicitly constructed or additional measures need to be taken for more complex differential operators
like operator-dependent coarsening or non-standard smoothers
\cite{esw05,tos01}. Algebraic multigrid methods are
often used as an alternative, in particular on unstructured meshes, but at a
somewhat higher cost in case more structure is available \cite{Gholami16}. In the
context of high-order methods, algebraic multigrid schemes need to be
carefully set up in order to not coarsen too aggressively due to the dense
coupling of the wider bases \cite{hmmo05}. Therefore, $p$-multigrid methods are
often considered in the high-order finite element context with a first transfer
to a low-order, usually linear, finite element basis \cite{lbl06} before
continuing with the algebraic hierarchy construction.

A second aspect that has not been covered by previous performance comparisons
is the fact that competitive high-order implementations in the primal
formulation are not based on matrices but rather most efficiently implemented
by matrix-free operator evaluation. A major reason for favoring matrix-free
methods is that most matrix-based iterative solvers are highly memory
bandwidth limited when executed on modern processors \cite{Patterson09} and
alternatives that access less memory can be faster also when performing
more arithmetic operations. In the low-order case with linear shape functions,
the most competitive matrix-free schemes typically rely on a stencil representation, e.g.~the
block structure in hierarchical hybrid grids \cite{Bergen06,Gmeiner15}. In
case of higher polynomial degrees, on-the-fly evaluation of cell and face
integrals with sum factorization is the preferred choice due to a low evaluation complexity.  Sum factorization is a technique established by the
spectral element community \cite{KS05,Kopriva09,O80} and used in
well-established codes like Nek5000 \cite{nek5000-web-page}, SPECFEM 3D
\cite{Komatitsch99}, or Nektar++ \cite{Vos10}. These methods are also popular
in the DG community \cite{Bastian2014,kk12}. At higher polynomial degrees $k$
in 3D, sum factorization for computing the integrals in matrix-vector products
of continuous finite elements has the same $\mathcal O(k^4)$ complexity per
element as the matrix after static condensation, but without a direct relation
to the stencil width of the matrix. This work uses the framework from
\cite{kk12} that has been demonstrated to perform particularly well yet being
flexible with respect to implementing generic differential operators and
providing equally optimized code paths for both continuous and discontinuous
Galerkin schemes. The implementation is available through the general-purpose
finite element library
\texttt{deal.II}\footnote{\texttt{http://www.dealii.org}, retrieved on
  November 13, 2016. The implementations used in this study are extensions of
  the step-37 and step-51 tutorial programs of \texttt{deal.II}, program URLs:
  \texttt{https://dealii.org/developer/doxygen/deal.II/step\_37.html} and
  \texttt{https://dealii.org/developer/doxygen/deal.II/step\_51.html}. A note
  to reviewers of the manuscript: The DG functionality has currently not yet
  been made available in the deal.II library but it will be merged during
  December 2016.  Thus, this comment will vanish in the final version of the
  manuscript.}  \cite{dealII84}. To the best of our knowledge, the numbers
presented in this study use the fasted CPU code of this kind available through a common
continuous and discontinuous finite element framework. The numbers are up to
an order of magnitude faster than what was reported in the recent study
\cite{Gholami16} with similar general-purpose geometric multigrid schemes,
reaching or even outperforming the specialized HPGMG \cite{hpgmg16} benchmark
code. For state-of-the-art implementations of distributed sparse matrix
algebra, the Trilinos\footnote{\texttt{http://www.trilinos.org}, retrieved on
  July 20, 2016} package \cite{trilinos} is used, providing a fair test bed
with mature implementations for both ends. The present study is novel in
comparing optimized sum factorization solvers to high-performance and
optimal complexity iterative solvers for HDG. Even though the authors in
\cite{Yakovlev16} (Remark 1) claim to use some form of sum factorization, the
implementation presented in this study leads to vastly different conclusions,
showing that previous work has missed some of the relevant aspects.

The remainder of this work is structured as follows. Section \ref{sec:discret}
introduces the Poisson equation and the discretizations with the continuous
finite element method, the symmetric interior penalty method, and the
hybridizable discontinuous Galerkin method. An analysis of matrix-vector
products of the various methods, including suggested alternatives for HDG as
opposed to the sparse trace matrix, are given in Section
\ref{sec:implementation}. The computational time to reach a certain level of
accuracy as well as performance metrics of the Poisson solvers are given in
Section \ref{sec:comparison}. Section \ref{sec:conclusions} summarizes our
findings.

\section{Discretization of Poisson's equation}
\label{sec:discret}

We consider the Poisson equation as a model problem for elliptic operators,
\begin{equation}\label{eq:poisson}
-\nabla \cdot (\kappa \nabla u) = f \quad \text{in } \Omega,
\end{equation}
where $u$ is the solution variable, $\kappa>0$ is a diffusion coefficient
bounded uniformly away from zero, and $\Omega$ is a bounded subset of
$d$-dimensional space $\mathbb{R}^d$. The domain boundary $\partial \Omega$ is
partitioned into a Dirichlet portion $\Gamma_\text{D}$ where $u = g_\text{D}$
and a Neumann portion $\Gamma_\text{N}$ where
$-\vec n \cdot \kappa \nabla u = g_\text{N}$ is prescribed,
respectively. Here, $\vec n$ denotes the unit outer normal vector on the
boundary $\Gamma_\text{N}$.

For discretization, we assume a tesselation $\mathcal T_h$ of the
computational domain $\Omega$ into $n_e$ elements $\Omega_{e}$, associated
with a mesh size parameter $h$. In this work, we assume a mesh consisting of
quadrilateral or hexahedral elements which allow for the most straight-forward
and efficient implementation of sum factorization, with lower proportionality
constants than tensorial techniques for triangles and tetrahedra
\cite{Vos10,schoeberl14}. All work known to the authors indicate that
matrix-based HDG schemes on quadrilaterals and hexahedra are at least as
efficient as on tetrahedra \cite{ksc12,ksmw15,Yakovlev16}, suggesting that our
results are unbiased in comparing against the fastest matrix-based options. We
assume an element $\Omega_e$ to be the image of the reference domain
$[-1,1]^d$ under a polynomial mapping of degree $l$, based on Gauss--Lobatto
support points that are placed according to a manifold description of the
computational domain. This enables high-order approximations of curved
boundaries and possibly also in the interior of $\Omega$. For the methods
described below, we denote the bilinear forms associated to integrals over the
elements of the triangulation as well as the faces by
\begin{equation}
  (a,b)_{\mathcal T_h} = \sum_{\Omega_e \in \mathcal T_h} \int_{\Omega_e} a\odot b\, d\vec x,\quad
  \left<a,b\right>_{\partial \mathcal T_h} = \sum_{\Omega_e \in \mathcal T_h}\sum_{F\in \text{faces}(\Omega_e)} \int_{F} a\odot b\, d\vec s,
\end{equation}
where $a,b$ can be scalar-valued, vector-valued, or tensor-valued quantities
and $\odot$ denotes the sum of the product in each component.

\subsection{Continuous Galerkin approximation}\label{sec:discret_fem}

We assume a polynomial approximation of the solution on elements from the space
\begin{equation}\label{eq:fem_space}
V_h^{\text{CG}} = \left\{v_h\in H^1(\Omega) : v_h|_{\Omega^e} \in \mathcal Q_k(\Omega^e)\ \forall \Omega_e \in \mathcal T_h\right\},
\end{equation}
where $\mathcal Q_k(\Omega_e)$ denotes the space of tensor product polynomials
of tensor degree $k$ on the element $\Omega_e$. In this work, we consider a
basis representation by Lagrange polynomials in the nodes of the $(k+1)$-point
Gauss--Lobatto--Legendre quadrature rule for well-conditioned high order
approximation \cite{KS05}. However, the exact form of the basis is immaterial,
as long as it is represented by a tensor product of 1D formulas. The solution
space is then restricted to the space $V_{h,{g_{\text{D}}}}^{\text{CG}}$ of
functions in $V_h^{\text{CG}}$ which satisfy the boundary condition
$g_{\text{D}}$ on $\Gamma_\text{D}$ by projection or interpolation.

The discrete finite element version of the Poisson equation \eqref{eq:poisson} is found by multiplication by a test function, integration over $\Omega$, integration by parts of the left hand side, and insertion of the Neumann boundary condition. The final weak form is to find a function $u_h \in V_{h,g_{\text{D}}}^{\text{CG}}$ such that
\begin{equation}\label{eq:poisson_fem}
\left(\nabla v_h, \kappa \nabla u_h\right)_{\mathcal T_h} = \left(v_h, f\right)_{\mathcal T_h} - \left<v_h, g_\text{N}\right>_{\mathcal \partial T_h \cap \Gamma_N}
\end{equation}
holds for all test functions $v_h \in V_{h,{0_{\text{D}}}}^{\text{CG}}$ that are
zero on the Dirichlet boundary.

On each element, the left-hand side gives rise to an element stiffness matrix
$K^e$ and the right-hand side to an element load vector $b^e$. These local quantities are assembled
into the global stiffness matrix $K$ and the load vector $b$ in the usual
finite element way, including the elimination of Dirichlet rows and columns. In
case of static condensation, the $(k-1)^d$ out of $(k+1)^d$ degrees of freedom
pertaining to basis functions with support on a single element are eliminated
by a Schur complement, reducing the final system size accordingly. We refer to
\cite{gmfh13,Yakovlev16} for details.

\subsection{Symmetric interior penalty discontinuous Galerkin discretization}\label{sec:discret_sip}

As a discontinuous Galerkin representative targeting the primal equation
amenable to sum factorization, we choose the symmetric interior penalty
(SIP) discontinuous Galerkin method \cite{abcm02}. In a discontinuous Galerkin
method, only $L_2$ regularity of the solution is required and no continuity
over element boundaries is enforced,
\begin{equation}\label{eq:dg_space}
V_h^{\text{DG}} = \left\{v_h\in L_2(\Omega) : v_h|_{\Omega_e} \in \mathcal Q_k(\Omega_e)\ \forall \Omega_e \in \mathcal T_h\right\}.
\end{equation}

On each element, the same steps as for continuous Galerkin in terms of
multiplication by test functions, integration over an element and integration
by parts, are taken. Due to the missing intra-element continuity, the terms
$-v_h \vec n \cdot \kappa \nabla u_h$ do not drop out over the interior faces
of the mesh and must be connected by a numerical flux on $\nabla u_h$.  The
first step is to take the average
$\frac 12\left(\nabla u_h^- + \nabla u_h^+\right)$ of the solution from both
elements $e^-$ and $e^+$ sharing a face. In order to ensure adjoint
consistency through a symmetric weak form, the term
$-\frac 12 (u_h^- \vec n^- + u_h^+ \vec n^+) \cdot \kappa \nabla v_h =- \frac
12 (u_h^- - u_h^+) \vec n^- \cdot \kappa \nabla v_h$
is added. This term is consistent with the
original equation because the difference $(u_h^- - u_h^+)$ is zero for the
analytic solution, see also \cite{Hesthaven08} for a derivation from a
first-order system. Finally, a penalty term
$\vec n v_h \kappa \sigma (u_h^--u_h^+)\vec n^-$ is added for ensuring
coercivity of discrete operator. The penalty parameter
$\sigma = (k+1)^2 \frac{d}{h}$ in $d$ dimensions depends on the inverse mesh
size $h$ on uniform meshes and is extended to general meshes by a formula
involving surface area and volume from \cite{Hillewaert13}. No tuning with
respect to $\sigma$ is done in this work. As documented in the literature
\cite{Hesthaven08}, condition numbers and multigrid performance would
deteriorate as $\sigma$ is increased.  This gives the following weak form for
the DG-SIP method,
\begin{equation}\label{eq:poisson_dgsip}
\begin{aligned}
&\left(\nabla v_h, \kappa \nabla u_h\right)_{\mathcal T_h}-\left<v_h \vec n, \kappa\frac{\nabla u_h^-+\nabla u_h^+}{2}\right>_{\partial \mathcal T_h} - \left<\frac{\nabla v_h}{2}, \kappa\vec n (u_h^--u_h^+)\right>_{\partial \mathcal T_h} \\
& \qquad \qquad \qquad + \left<v_h, \kappa\sigma (u_h^--u_h^+)\right>_{\partial \mathcal T_h} = \left(v_h, f\right)_{\mathcal T_h}
\end{aligned}
\end{equation}
which is to hold for all test functions $v_h$ in the space
$V_h^{\text{DG}}$. Note that the bilinear forms
$\left<\cdot,\cdot\right>_{\partial \mathcal T_h}$ visit each interior face
twice with opposite directions of the normal vector $\vec n$, resulting in a symmetric
weak form. Boundary conditions are imposed by defining
suitable extension values $u^+$ in terms of the boundary condition and the
inner solution value $u^-$,
\begin{equation}\label{eq:bc_dgsip}
\begin{aligned}
&u^+ = -u^- + 2 g_\text{D}, &&\nabla u^+ = \nabla u^-, && \text{on Dirichlet boundaries,}\\
&u^+ = u^-, &&\nabla u^+\cdot \vec n = -\nabla u^-\cdot \vec n  - 2 \frac{g_\text{N}}{\kappa}, && \text{on Neumann boundaries.}
\end{aligned}
\end{equation}
Thus, additional contributions of known quantities arise that are eventually
moved to the right-hand side of the final linear system.

\subsection{Hybridizable discontinuous Galerkin discretization}\label{sec:discret_hdg}

For the hybri\-dizable discontinuous Galerkin (HDG) discretization \cite{cgl09}, the
Poisson equation \eqref{eq:poisson} is rewritten as a first-order system by
introducing a flux variable $\vec q = -\kappa \nabla u$ in the
equation $\nabla \cdot \vec q = f$. The discrete solution spaces are
$V_h^{\text{DG}}$ for $u_h$ and $\left(V_h^{\text{DG}}\right)^d$ for
$\vec q_h$. An additional trace variable $\lambda_h$ that approximates $u_h$
on the interface between elements is introduced. It is defined by polynomials on the mesh skeleton
\begin{equation}\label{eq:trace_space}
M_h^{\text{tHDG}} = \left\{\mu_h\in L_2\left(\mathcal F_h\right) : \mu_h|_{F} \in \mathcal Q_k(F)\ \forall \text{ faces }F \in \mathcal F_h\right\},
\end{equation}
where $\mathcal F_h$ denotes the collection of all faces in the discretization
$\mathcal T_h$. The functions in $M_h^\text{tHDG}$ are discontinuous between
faces (i.e., over vertices in 2D, over vertices and edges in 3D). Besides
tensor product polynomials $\mathcal Q_k(F)$, we will also consider
polynomials of complete degree $k$, $\mathcal P_k(F)$, in combination with the
space $\mathcal P_k(\Omega_e)$, easily permitted by the discontinuous
formulation on $F$ and $\Omega_e$.

For deriving the HDG weak form, the system is multiplied by test functions
$\vec w_h$, $v_h$, integrated over element $\Omega_e$, and gradient and
divergence terms are integrated by parts. For the numerical fluxes, we add a
new variable for the first flux, $\hat{u} = \lambda_h$, while the second flux
is set to $\hat{\vec q} = \vec q_h + \tau (u_h - \lambda_h)\vec n$.  The
system is closed by enforcing continuity on the second numerical flux
$\hat{\vec q}$, which ensures conservativity. The final weak form for the HDG method is to find the values
$\vec q_h \in (V_h^{\text{DG}})^d$,  $u_h\in V_h^{\text{DG}}$, and $\lambda_h\in M_{h,g_{\text{D}}}^\text{tHDG}$, such that
\begin{equation}\label{eq:hdg}
\begin{aligned}
\left(\vec w_h, \kappa^{-1} \vec q_h\right)_{\mathcal T_h} - \left(\nabla \cdot \vec w_h, u_h\right)_{\mathcal T_h} + \left<\vec w_h, \vec n \lambda_h\right>_{\partial \mathcal T_h} &= 0,\\
\left(\nabla v_h, \vec q_h\right)_{\mathcal T_h}  + \left<v_h, \vec n\cdot \vec q_h + \tau(u_h - \lambda_h)\right>_{\partial \mathcal T_h} &= -\left(v_h, f\right)_{\mathcal T_h},\\
-\left<\mu_h, \vec n\cdot \vec q_h + \tau(u_h - \lambda_h)\right>_{\partial \mathcal T_h} &= -\left<\mu_h, g_\text{N}\right>_{\partial \mathcal T_h \cap \Gamma_\text{N}},
\end{aligned}
\end{equation}
holds for all test functions $\vec w_h\in (V_h^{\text{DG}})^d$,
$v_h\in V_h^{\text{DG}}$, and $\mu_h \in M_{h,0}^\text{tHDG}$.
Dirichlet conditions are imposed strongly on the trace space $M_h^\text{tHDG}$
by projection.

The parameter $\tau$ ensures stability if chosen
$\tau > 0$ \cite{npc09}. However, selecting $\tau \propto \kappa$ is advantageous
because it gives optimal convergence rates $k+1$ in both the solution and the
flux. An element-by-element post-processing can then be performed to recover a
solution $u^*$ that converges at rate $k+2$ \cite{cgl09}. Even though
superconvergence has not been proved for quadrilaterals and hexahedra (which
need special projection properties, as opposed to tetrahedra
\cite{Cockburn09}), we observed superconvergence for all constant-coefficient
elliptic test cases with shape-regular but otherwise arbitrary quadrilateral
and hexahedral meshes, including the test case from \cite{ksc12} where the
authors report only rates $k+1$ on quadrilaterals. Note that the
superconvergence of HDG essentially contributes to the high efficiency as
compared to other methods documented in previous studies
\cite{harp13,ksc12,Yakovlev16}.

Following the notation in
\cite{npc09}, the individual terms in Equation \eqref{eq:hdg} are expanded in
terms of the basis functions and put into matrix-vector form. The system
reads
\begin{equation}\label{eq:hdg_matrix}
\begin{pmatrix}A & B^\mathrm T & C^\mathrm T \\ B & D & G^\mathrm T \\ -C & -G & H\end{pmatrix}\begin{pmatrix}\vec Q_h \\ \vec U_h \\ \Lambda_h \end{pmatrix} = \begin{pmatrix} \vec 0 \\ \vec R_f \\ \vec R_{g_\text{N}}\end{pmatrix}.
\end{equation}
The upper left $2\times 2$ block is block-diagonal over elements and can be
condensed out before solving the linear system. Thus, only a symmetric
positive definite linear system $ K \Lambda_h = \vec R $  needs
to be solved with
\begin{equation}\label{eq:hdg_matrix_condensed}
 K = H + \big( C \  G\big) \begin{pmatrix} A & B^\mathrm T \\ B & D \end{pmatrix}^{-1} \begin{pmatrix} C^\mathrm T \\ G^\mathrm T \end{pmatrix},\ \vec R = \vec R_{g_\text{N}} + \big( C \  G \big) \begin{pmatrix} A & B^\mathrm T \\ B & D \end{pmatrix}^{-1} \begin{pmatrix}\vec 0 \\ \vec R_f\end{pmatrix}.
\end{equation}

\section{Performance of operator evaluation}\label{sec:implementation}

In this work, we consider high-performance matrix-free evaluation of
matrix-vector products whenever possible, relying on fast
integration facilities established in spectral elements
\cite{KS05,Kopriva09}. While these methods had originally only been used in
the high-degree context with $k\geq 4$, recent high-performance realizations
taking the architecture of modern parallel computers into account have shown
that matrix-free methods outperform sparse matrix kernels by several times already for $k=2$ on
quadrilaterals and hexahedra \cite{brown10,kk12}.

In a matrix-free setting, the matrix-vector product is interpreted as a weak
form that is tested by all basis functions in an element-by-element way. This
corresponds to computing a residual vector on each element that is assembled
into the global solution vector, given a function $u_h$ associated to the
input vector \cite{kk12}. If denote by $\vec z = L \vec y$ the global operator
evaluation and by $\vec z_e = L_e \vec y_e$ the contribution on element
$\Omega_e$, the $i$-th component of the matrix-vector product is given by
quadrature
\begin{equation}\label{eq:mf_eval}
  (z_e)_i = \int_{\Omega_e} \kappa \nabla \phi_i \nabla u_h^{y_e} d\vec x \approx \sum_{q}  \nabla_{\vec \xi} \phi_i \left(J^e\right)^{-1}\left(w_q \text{det}(J^e) \kappa\right) \left(J^e\right)^{-\mathrm T} \nabla_{\vec \xi} u_h^{y_e}.
\end{equation}
In this expression, $u_h^{y_e}$ is the finite element function associated to
the nodal solution values $\vec y_e$, $J$ is the Jacobian of the
transformation from the reference to the real cell, and $w_q$ the quadrature
weight. In this work, we use Gauss--Legendre quadrature with $k+1$ points per
coordinate direction which exactly evaluates integrals on Cartesian
geometries. The error on curved meshes does not affect convergence orders for
elliptic problems in the current setting apart from the usual variational
crime \cite{brenner08}.

For the interpolation of the nodal values $y_e$ to the quadrature points for
$\nabla_{\vec \xi} u_h^{y_e}$ in Equation \eqref{eq:mf_eval} as well as the
multiplication by the test function gradients $\nabla_{\vec \xi} \phi_i$ for
all test functions $i=1,\ldots,(k+1)^d$, the sum factorization technique is
used \cite{kk12}. This approach replaces the direct interpolation over all
points in $d$ dimensions by a series of $d$ one-dimensional
interpolations. This reduces the evaluation complexity per element from
$\mathcal O((k+1)^{2d})$ operations in the naive matrix-vector product with a
Kronecker matrix to $\mathcal O(d(k+1)^{d+1})$ operations. The complexity of
DG face integrals is only $\mathcal O((k+1)^{d})$, i.e., linear in the number
of unknowns \cite{Bastian2014}.

Our implementation \cite{KK11,kk12} is available through the \texttt{deal.II}
finite element library \cite{dealII84} and specifically targets modern
computer architecture where access to main memory is usually the bottleneck
for PDE-based kernels. This means that all operations on an element according
to Equation \eqref{eq:mf_eval} are done in close temporal proximity to service
the sum factorization kernels from fast L1 caches. The complexity
of one-dimensional kernels in sum factorization is further cut down into half by a
high-degree optimization based on the even-odd decomposition
\cite{Kopriva09}. A memory optimization is applied for Cartesian and affine
meshes where the Jacobian of the transformation from the reference to the real
cell is constant throughout the whole cell and needs only be kept once. In
case cells are not affine, it is fastest \cite{kk12} to pre-compute the
Jacobian on all quadrature points and all cells prior to solving linear
systems and loading the coefficients in each matrix-vector product, despite
the relatively high memory transfer. The experiments below consider both the
memory-intensive general mesh case and the simple Cartesian mesh case as they
show different performance. For sparse linear algebra, the Epetra backend of
Trilinos is used due to its mature state \cite{trilinos}.

\subsection{Variants of matrix-vector products}

Since the continuous finite element method, DG-SIP, and the HDG trace
matrix all involve different matrix sizes but the error is most closely
related to the number of elements, the most appropriate metric for comparison
would be the cost per element. On the other hand, we want to acknowledge the
ability of higher order basis functions to use coarser meshes. Thus, we
report the numbers in this section as the cost per element divided by $k^d$,
the number of unique degrees of freedom (DoF) per element on a continuous finite
element space. In other words, if a matrix-vector product on $n_e$ elements of
degree $k$ takes $t_\text{mv}$ seconds, we report the quantity
\begin{equation}\label{eq:equiv_dofs}
\text{Equivalent DoFs/s} = \frac{n_e k^d}{t_\text{mv}},
\end{equation}
uniformly across all discretization schemes.
Obviously, this selection does not realistically represent discretization accuracy,
a topic we postpone to Section \ref{sec:comparison}. The measurements in this section have been performed on
a fully utilized node of dual-socket Intel Xeon E5-2690 v4 (Broadwell)
processors with $2\times 14$ cores running at 2.6 GHz and eight memory
channels. By using the full node, a fair balance between arithmetic bound
kernels (sum factorization) and memory bound kernels (sparse matrix-vector
products) is achieved. Reported memory bandwidth of this setup is
approximately 130~GB/s in the STREAM triad benchmark \cite{MC07} or sparse
matrix-vector products, whereas the theoretical arithmetic peak is 940~GFLOP/s
when measured at the AVX base frequency of 2.1~GHz. All C++ code has been
compiled with the GNU compiler \texttt{gcc}, version 6.1, with optimization
target AVX2 (Haswell). The minimum runtime out of five experiments is presented.

Fig.~\ref{fig:matvec} compares the three discretization methods and several
implementations of the matrix-vector product:
\begin{itemize}
\item CG matrix-free: A standard continuous Galerkin approximation
  with tensor product basis functions of degree $k$ and Gaussian quadrature on
  $(k+1)^d$ points according to Sec.~\ref{sec:discret_fem} evaluated in a matrix-free way.
\item CG static condensation matrix: Uses the most efficient sparse matrix
  representation obtained by static condensation of the $(k-1)^d$ cell-interior degrees
  of freedom in continuous elements \cite{harp13}.
\item DG-SIP matrix-free: Symmetric interior penalty DG method according to
  Sec.~\ref{sec:discret_sip} with matrix-free implementation through sum
  factorization. No matrix is considered due to its low efficiency \cite{harp13}.
\item HDG trace matrix: Sparse matrix-vector product with the trace matrix.
\item HDG trace matrix post: Takes the increased accuracy of HDG by
  super-convergent post-processing into account. For this label, the data from ``HDG trace
  matrix'' measured at degree $k-1$ is reported in terms of the equivalent
  degrees of freedom in Equation~\eqref{eq:equiv_dofs} of one degree higher, $k^d$.
\item HDG trace matrix-free: This approach considers an alternative
  implementation of the trace matrix-vector product where the matrix system
  \eqref{eq:hdg_matrix_condensed} is expanded in terms of all contributing
  matrices rather than explicitly forming the Schur complement. This allows
  most matrix-vector products to be implemented in a matrix-free way, using a
  scheme originally proposed in \cite{ksmw15} for fast computation of HDG
  residuals:
  \begin{enumerate}
  \item Matrix-free multiplication by
    $\begin{pmatrix} C^\mathrm T \\ G^\mathrm T\end{pmatrix}$ on input
    $\Lambda_h$.
  \item Application of the inverse matrix
    $\begin{pmatrix}A & B^\mathrm T \\ B & D\end{pmatrix}^{-1}$ by an inner
    Schur complement that computes the nodal values $\vec U_h$, $\vec Q_h$
    given the vectors $\vec R_q$, $\vec R_u$ from step 1:
    \begin{align*}
        \vec U_h &= \left(D - B A^{-1} B^\mathrm T \right)^{-1} \left( \vec R_u - B {A}^{-1} \vec{R}_q  \right), \\
        \vec Q_h &= A^{-1} \left( \vec{R}_q - B^\mathrm T \vec U_h \right).
    \end{align*}
    In this equation, the matrix-vector multiplications by $A^{-1}, B$, and
    $B^\mathrm T$ can all be implemented by matrix-free evaluation, including
    the inverse vector mass matrix $A^{-1}$ for which fast sum factorization
    techniques exist \cite{ksmw15}. Only the $(k+1)^d\times (k+1)^d$ matrix
    $(D - B A^{-1} B^\mathrm T)^{-1}$ needs to be explicitly stored.
    \item Matrix-free multiplication by
      $\begin{pmatrix}C & G\end{pmatrix}$ and the mass matrix $H$.
  \end{enumerate}
  This approach is counter-intuitive because it defeats the original purpose
  of the static condensation-type elimination of degrees of freedom.  This
  approach is slower than the trace matrix-vector products for all degrees in
  2D according to Fig.~\ref{fig:matvec}. However, up to twice as high
  performance until degree five is recorded in 3D due to much reduced memory
  transfer. The higher complexity $(k+1)^{6}$ per element of the
  multiplication by the dense matrix $(D - B A^{-1} B^\mathrm T)^{-1}$
  dominates over the trace complexity $(k+1)^4$ for $k>5$.
\item HDG mixed matrix-free: This approach replaces the condensation into the
  trace matrix by applying the numerical fluxes in a more classical DG way in
  terms of a system in mixed form in degrees of freedom for $u$ and $\vec q$,
  \begin{equation*}
  \begin{pmatrix}A & B^\mathrm T \\ B & D\end{pmatrix} + \begin{pmatrix}C^\mathrm T \\ G^\mathrm T \end{pmatrix} H^{-1} \begin{pmatrix} C & G \end{pmatrix},
  \end{equation*}
  This is the form taken for explicit time integration in many DG schemes,
  including HDG \cite{ksmw15}. All matrix-vector products can be performed by
  sum factorization, including evaluation of $H^{-1}$ that is an inverse face
  mass matrix \cite{ksmw15} or alternatively can be implemented by point-wise
  fluxes.
\end{itemize}

\subsection{Throughput analysis}

The throughput of the matrix-vector products in terms of the number of
equivalent degrees of freedom processed per second is measured at large matrix
sizes of around 10 million where the solution vectors (and all other global
data) need to be fetched from main memory rather than caches.

\begin{figure}
  \centering
  \begin{tikzpicture}
    \begin{semilogyaxis}[
      width=0.52\textwidth,
      height=0.45\textwidth,
      title style={at={(1,0.948)},anchor=north east,draw=black,fill=white,font=\footnotesize},
      title={2D Cartesian mesh},
      ylabel={Equivalent DoFs/s},
      y label style={at={(0.1,0.5)}},
      legend columns = 3,
      legend to name=legendMV2d,
      ymin=4e7,
      ymax=5e9,
      xmin=1,
      xmax=8,
      xtick={1,2,3,4,5,6,7,8},
      tick label style={font=\scriptsize},
      label style={font=\scriptsize},
      legend style={font=\scriptsize},
      grid
      ]
      \addplot table [x index=0, y expr={\thisrowno{1}*\thisrowno{0}^2/\thisrowno{3}}] {
  degree  cells   dofs    matvec
  1     4194304  4198401  5.5294e-03
  2     4194304  16785409 8.8915e-03
  3     1048576  9443329  4.1563e-03
  4     1048576  16785409 6.8622e-03
  5     1048576  26224641 1.0493e-02
  6     1048576  37761025 1.4608e-02
  7     1048576  51394561 2.0382e-02
  8     1048576  67125249 2.6955e-02
      };
\addlegendentry{CG matrix-free};
\pgfplotsset{cycle list shift=3}
      \addplot table [x index=0, y expr={\thisrowno{1}*\thisrowno{0}^2/\thisrowno{3}}] {
  degree  cells   dofs    matvec
  1     4194304  4198401  5.8559e-03
  2     4194304  16785409 2.6174e-02
  3     1048576  9443329  1.4483e-02
  4     1048576  16785409 2.5456e-02
  5     1048576  26224641 3.9382e-02
  6     1048576  37761025 5.6013e-02
  7     1048576  51394561 7.6411e-02
  8     1048576  67125249 9.8695e-02
      };
\addlegendentry{CG stat.~cond.~matrix};
\pgfplotsset{cycle list shift=-1}
      \addplot table [x index=0, y expr={\thisrowno{1}*\thisrowno{0}^2/\thisrowno{3}}] {
  degree  cells   dofs    matvec
  1     4194304  16777216 2.5078e-02
  2     4194304  37748736 3.8027e-02
  3     1048576  16777216 1.3741e-02
  4     1048576  26214400 1.9627e-02
  5     1048576  37748736 2.6701e-02
  6     1048576  51380224 3.5883e-02
  7     1048576  67108864 4.6258e-02
  8     1048576  84934656 5.8920e-02
      };
\addlegendentry{DG-SIP matrix-free};
\pgfplotsset{cycle list shift=0}
      \addplot table [x index=0, y expr={\thisrowno{1}*\thisrowno{0}^2/\thisrowno{3}}] {
  degree  cells   dofs    matvec
  1     4194304  16785408 3.1900e-02
  2     4194304  25178112 6.5997e-02
  3     1048576  8396800  2.7699e-02
  4     1048576  10496000 4.2058e-02
  5     1048576  12595200 5.9246e-02
  6     1048576  14694400 7.9349e-02
  7     1048576  16793600 1.0272e-01
  8     1048576  18892800 1.2900e-01
      };
\addlegendentry{HDG trace matrix};
      \addplot[mark=triangle*,dashed,black,every mark/.append style={solid}] table [x expr={\thisrowno{0}+1}, y expr={\thisrowno{1}*(\thisrowno{0}+1)^2/\thisrowno{3}}] {
  degree  cells   dofs    matvec
  1     4194304  16785408 3.1900e-02
  2     4194304  25178112 6.5997e-02
  3     1048576  8396800  2.7699e-02
  4     1048576  10496000 4.2058e-02
  5     1048576  12595200 5.9246e-02
  6     1048576  14694400 7.9349e-02
  7     1048576  16793600 1.0272e-01
  8     1048576  18892800 1.2900e-01
      };
\addlegendentry{HDG trace matrix post};
      \addplot[mark=halfcircle,densely dashdotted,every mark/.append style={solid}] table [x index=0, y expr={\thisrowno{1}*\thisrowno{0}^2/\thisrowno{3}}] {
  degree  cells   dofs     matvec
  1     4194304  16785408  5.4171e-02
  2     4194304  25178112  1.0116e-01
  3     1048576  8396800   4.8917e-02
  4     1048576  10496000  8.8939e-02
  5     1048576  12595200  1.5049e-01
  6     1048576  14694400  2.4529e-01
  7     1048576  16793600  3.8736e-01
  8     1048576  18892800  5.8872e-01
      };
\addlegendentry{HDG trace matrix-free};
      \addplot[mark=oplus*,densely dotted,mark options={fill=black!40},every mark/.append style={solid}] table [x index=0, y expr={\thisrowno{1}*\thisrowno{0}^2/\thisrowno{3}}] {
  degree  cells   dofs     matvec
  1     4194304  50331648  4.1420e-02
  2     4194304  113246208 7.3878e-02
  3     1048576  50331648  3.2274e-02
  4     1048576  78643200  4.8264e-02
  5     1048576  113246208 6.7956e-02
  6     1048576  154140672 9.0828e-02
  7     1048576  201326592 1.1795e-01
  8     1048576  254803968 1.4876e-01
      };
\addlegendentry{HDG mixed matrix-free};
    \end{semilogyaxis}
\end{tikzpicture}
\quad
  \begin{tikzpicture}
    \begin{semilogyaxis}[
      width=0.52\textwidth,
      height=0.45\textwidth,
      title style={at={(1,0.948)},anchor=north east,draw=black,fill=white,font=\footnotesize},
      title={2D curved mesh},
      ymin=4e7,
      ymax=5e9,
      xmin=1,
      xmax=8,
      xtick={1,2,3,4,5,6,7,8},
      tick label style={font=\scriptsize},
      label style={font=\scriptsize},
      grid
      ]
      \addplot table [x index=0, y expr={\thisrowno{1}*\thisrowno{0}^2/\thisrowno{3}}] {
  degree  cells   dofs    matvec
  1     12582912 12595200 2.4355e-02
  2     3145728  12595200 1.4995e-02
  3     3145728  28329984 2.8375e-02
  4     786432   12595200 1.1354e-02
  5     786432   19676160 1.6918e-02
  6     786432   28329984 2.3409e-02
  7     786432   38556672 3.1264e-02
  8     786432   50356224 3.9700e-02
      };
\pgfplotsset{cycle list shift=3}
      \addplot table [x index=0, y expr={\thisrowno{1}*\thisrowno{0}^2/\thisrowno{3}}] {
  degree  cells   dofs    matvec
  1     12582912 12595200 1.7889e-02
  2     3145728  12595200 1.9849e-02
  3     3145728  28329984 4.3928e-02
  4     786432   12595200 1.9627e-02
  5     786432   19676160 3.0332e-02
  6     786432   28329984 4.3347e-02
  7     786432   38556672 5.8433e-02
  8     786432   50356224 7.5654e-02
      };
\pgfplotsset{cycle list shift=-1}
      \addplot table [x index=0, y expr={\thisrowno{1}*\thisrowno{0}^2/\thisrowno{3}}] {
  degree  cells   dofs    matvec
  1     3145728  12582912 2.3524e-02
  2     3145728  28311552 3.9629e-02
  3     3145728  50331648 6.1435e-02
  4     786432   19660800 2.2498e-02
  5     786432   28311552 3.1049e-02
  6     786432   38535168 4.1272e-02
  7     786432   50331648 5.3212e-02
  8     786432   63700992 6.6854e-02
      };
\pgfplotsset{cycle list shift=0}
      \addplot table [x index=0, y expr={\thisrowno{1}*\thisrowno{0}^2/\thisrowno{3}}] {
  degree  cells   dofs    matvec
  1     3145728  12595200 2.4076e-02
  2     3145728  18892800 5.0137e-02
  3     3145728  25190400 8.4989e-02
  4     786432    7879680 3.2165e-02
  5     786432    9455616 4.5245e-02
  6     786432   11031552 6.1283e-02
  7     786432   12607488 7.9083e-02
  8     786432   14183424 9.9639e-02
      };
      \addplot[mark=triangle*,dashed,black,every mark/.append style={solid}] table [x expr={\thisrowno{0}+1}, y expr={\thisrowno{1}*(\thisrowno{0}+1)^2/\thisrowno{3}}] {
  degree  cells   dofs    matvec
  1     3145728  12595200 2.4076e-02
  2     3145728  18892800 5.0137e-02
  3     3145728  25190400 8.4989e-02
  4     786432    7879680 3.2165e-02
  5     786432    9455616 4.5245e-02
  6     786432   11031552 6.1283e-02
  7     786432   12607488 7.9083e-02
  8     786432   14183424 9.9639e-02
      };
      \addplot[mark=halfcircle,densely dashdotted,every mark/.append style={solid}] table [x index=0, y expr={\thisrowno{1}*\thisrowno{0}^2/\thisrowno{3}}] {
  degree  cells   dofs     matvec
  1     3145728  12595200  4.3233e-02
  2     3145728  18892800  8.5575e-02
  3     3145728  25190400  1.6453e-01
  4     786432   7879680   7.2311e-02
  5     786432   9455616   1.2374e-01
  6     786432   11031552  2.0031e-01
  7     786432   12607488  3.0982e-01
  8     786432   14183424  4.6509e-01
      };
      \addplot[mark=oplus*,densely dotted,mark options={fill=black!40},every mark/.append style={solid}] table [x index=0, y expr={\thisrowno{1}*\thisrowno{0}^2/\thisrowno{3}}] {
  degree  cells   dofs     matvec
  1     3145728  37748736  3.5019e-02
  2     3145728  84934656  6.4896e-02
  3     3145728 150994944  1.1184e-01
  4     786432   58982400  4.2943e-02
  5     786432   84934656  6.0660e-02
  6     786432  115605504  8.0728e-02
  7     786432  150994944  1.0475e-01
  8     786432  191102976  1.3204e-01
      };
    \end{semilogyaxis}
\end{tikzpicture}
\\
  \begin{tikzpicture}
    \begin{semilogyaxis}[
      width=0.52\textwidth,
      height=0.45\textwidth,
      title style={at={(1,0.948)},anchor=north east,draw=black,fill=white,font=\footnotesize},
      title={3D Cartesian mesh},
      xlabel={Polynomial degree},
      ylabel={Equivalent DoFs/s},
      x label style={at={(0.5,0.04)}},
      y label style={at={(0.1,0.5)}},
      ymin=1e7,
      ymax=3e9,
      xmin=1,
      xmax=8,
      xtick={1,2,3,4,5,6,7,8},
      tick label style={font=\scriptsize},
      label style={font=\scriptsize},
      grid
      ]
      \addplot table [x index=0, y expr={\thisrowno{1}*\thisrowno{0}^3/\thisrowno{3}}] {
  degree  cells   dofs    matvec
  1     11239424 11390625 2.7667e-02
  2     2097152  16974593 1.4309e-02
  3     512000   13997521 1.0738e-02
  4     262144   16974593 1.2101e-02
  5     110592   13997521 1.0354e-02
  6     32768    7189057  5.4458e-03
  7     32768    11390625 8.9623e-03
  8     21952    11390625 9.4777e-03
      };
\pgfplotsset{cycle list shift=3}
      \addplot table [x index=0, y expr={\thisrowno{1}*\thisrowno{0}^3/\thisrowno{3}}] {
  degree  cells   dofs    matvec
  1     11239424 11390625 3.7634e-02
  2     2097152  16974593 1.0882e-01
  3     512000   13997521 1.3502e-01
  4     262144   16974593 2.1315e-01
  5     110592   13997521 2.1951e-01
  6     32768    7189057  1.3700e-01
  7     32768    11390625 2.5105e-01
  8     21952    11390625 2.8383e-01
      };
\pgfplotsset{cycle list shift=-1}
      \addplot table [x index=0, y expr={\thisrowno{1}*\thisrowno{0}^3/\thisrowno{3}}] {
  degree  cells   dofs    matvec
  1     11239424 89915392 1.7608e-01
  2     2097152  56623104 7.1049e-02
  3     262144   16777216 2.0730e-02
  4     262144   32768000 4.0241e-02
  5     110592   23887872 3.1817e-02
  6     32768    11239424 1.6657e-02
  7     32768    16777216 2.5784e-02
  8     21952    16003008 2.7062e-02
      };
\pgfplotsset{cycle list shift=0}
      \addplot table [x index=0, y expr={\thisrowno{1}*\thisrowno{0}^3/\thisrowno{3}}] {
  degree  cells   dofs    matvec
  1    11239424 135475200 6.7959e-01
  2     2097152  57065472 6.0332e-01
  3     262144   12779520 2.3427e-01
  4     262144   19968000 5.5929e-01
  5     110592   12192768 4.8504e-01
  6     32768     4967424 2.6433e-01
  7     32768     6488064 4.4736e-01
  8     21952     5524848 4.7292e-01
      };
      \addplot[mark=triangle*,dashed,black,every mark/.append style={solid}] table [x expr={\thisrowno{0}+1}, y expr={\thisrowno{1}*(\thisrowno{0}+1)^3/\thisrowno{3}}] {
  degree  cells   dofs    matvec
  1     2097152  25362432 1.2796e-01
  2     884736   24136704 2.5740e-01
  3     262144   12779520 2.3427e-01
  4     262144   19968000 5.5929e-01
  5     110592   12192768 4.8504e-01
  6     32768     4967424 2.6433e-01
  7     32768     6488064 4.4736e-01
  8     13824     3499200 3.0505e-01
      };
      \addplot[mark=halfcircle,densely dashdotted,every mark/.append style={solid}] table [x index=0, y expr={\thisrowno{1}*\thisrowno{0}^3/\thisrowno{3}}] {
  degree  cells   dofs     matvec
  1     2097152  25362432  8.6476e-02
  2     884736   24136704  1.2496e-01
  3     262144   12779520  1.2478e-01
  4     262144   19968000  3.7289e-01
  5     110592   12192768  4.2075e-01
  6     32768     4967424  3.1717e-01
  7     32768     6488064  6.7836e-01
  8     13824     3499200  5.5947e-01
      };
      \addplot[mark=oplus*,densely dotted,mark options={fill=black!40},every mark/.append style={solid}] table [x index=0, y expr={\thisrowno{1}*\thisrowno{0}^3/\thisrowno{3}}] {
  degree  cells   dofs     matvec
  1     2097152  67108864  7.2134e-02
  2     884736   95551488  8.7822e-02
  3     262144   67108864  6.4485e-02
  4     262144  131072000  1.2252e-01
  5     110592   95551488  9.7681e-02
  6     32768    44957696  5.3768e-02
  7     32768    67108864  8.2439e-02
  8     21952    64012032  8.5974e-02
      };
    \end{semilogyaxis}
\end{tikzpicture}
\quad
  \begin{tikzpicture}
    \begin{semilogyaxis}[
      title style={at={(1,0.948)},anchor=north east,draw=black,fill=white,font=\footnotesize},
      width=0.52\textwidth,
      height=0.45\textwidth,
      title={3D curved mesh},
      xlabel={Polynomial degree},
      x label style={at={(0.5,0.05)}},
      ymin=1e7,
      ymax=3e9,
      xmin=1,
      xmax=8,
      xtick={1,2,3,4,5,6,7,8},
      tick label style={font=\scriptsize},
      label style={font=\scriptsize},
      grid
      ]
      \addplot table [x index=0, y expr={\thisrowno{1}*\thisrowno{0}^3/\thisrowno{3}}] {
  degree  cells   dofs    matvec
  1     3145728 3195010  2.0744e-02
  2     3145728 25362690 7.4294e-02
  3     393216  10727618 2.4371e-02
  4     393216  25362690 5.5412e-02
  5     49152   6220962  1.1550e-02
  6     49152   10727618 1.8706e-02
  7     49152   17009890 2.8109e-02
  8     6144    3195010  5.7757e-03
      };
\pgfplotsset{cycle list shift=3}
      \addplot table [x index=0, y expr={\thisrowno{1}*\thisrowno{0}^3/\thisrowno{3}}] {
  degree  cells   dofs    matvec
  1     3145728 3195010  1.0717e-02
  2     3145728 25362690 1.6436e-01
  3     393216  10727618 1.0331e-01
  4     393216  25362690 3.2509e-01
  5     49152   6220962  9.9635e-02
  6     49152   10727618 2.0565e-01
  7     49152   17009890 3.8032e-01
  8     6144    3195010  7.9736e-02
      };
\pgfplotsset{cycle list shift=-1}
      \addplot table [x index=0, y expr={\thisrowno{1}*\thisrowno{0}^3/\thisrowno{3}}] {
  degree  cells   dofs    matvec
  1     3145728 25165824 7.2502e-02
  2     3145728 84934656 1.8667e-01
  3     393216  25165824 5.5024e-02
  4     393216  49152000 1.0439e-01
  5     49152   10616832 2.4459e-02
  6     49152   16859136 3.8663e-02
  7     49152   25165824 5.7023e-02
  8     6144    4478976  1.2699e-02
      };
\pgfplotsset{cycle list shift=0}
      \addplot table [x index=0, y expr={\thisrowno{1}*\thisrowno{0}^3/\thisrowno{3}}] {
  degree  cells   dofs    matvec
  1     3145728 37945344 1.9802e-01
  2     3145728 85377024 9.1729e-01
  3     393216  19070976 3.5740e-01
  4     393216  29798400 8.4468e-01
  5     49152   5419008  2.1872e-01
  6     49152   7375872  4.0089e-01
  7     49152   9633792  6.7712e-01
  8     6144    1555200  1.3259e-01
      };
      \addplot[mark=triangle*,dashed,black,every mark/.append style={solid}] table [x expr={\thisrowno{0}+1}, y expr={\thisrowno{1}*(\thisrowno{0}+1)^3/\thisrowno{3}}] {
  degree  cells   dofs    matvec
  1     3145728 37945344 1.9802e-01
  2     3145728 85377024 9.1729e-01
  3     393216  19070976 3.5740e-01
  4     393216  29798400 8.4468e-01
  5     49152   5419008  2.1872e-01
  6     49152   7375872  4.0089e-01
  7     49152   9633792  6.7712e-01
  8     6144    1555200  1.3259e-01
      };
      \addplot[mark=halfcircle,densely dashdotted,every mark/.append style={solid}] table [x index=0, y expr={\thisrowno{1}*\thisrowno{0}^3/\thisrowno{3}}] {
  degree  cells   dofs     matvec
  1     3145728 37945344 1.5234e-01
  2     3145728 85377024 5.1668e-01
  3     393216  19070976 2.1354e-01
  4     393216  29798400 6.0683e-01
  5     49152   5419008  2.0428e-01
  6     49152   7375872  5.0209e-01
  7     49152   9633792  1.0549e+00
  8     6144    1555200  2.6140e-01
      };
      \addplot[mark=oplus*,densely dotted,mark options={fill=black!40},every mark/.append style={solid}] table [x index=0, y expr={\thisrowno{1}*\thisrowno{0}^3/\thisrowno{3}}] {
  degree  cells   dofs     matvec
  1     3145728 100663296 1.2642e-01
  2     3145728 339738624 3.6195e-01
  3     393216  100663296 1.1639e-01
  4     393216  196608000 2.2239e-01
  5     49152   42467328  5.9396e-02
  6     49152   67436544  9.4426e-02
  7     49152   100663296 1.4534e-01
  8     6144    17915904  3.5387e-02
      };
    \end{semilogyaxis}
\end{tikzpicture}
\\
\ref{legendMV2d}
\caption{Number of degrees of freedom processed per second on 28 Broadwell cores
  as a function of the polynomial degree for various
  discretizations and implementations of the matrix-vector
  product of the 2D and 3D Laplacian.}
\label{fig:matvec}
\end{figure}
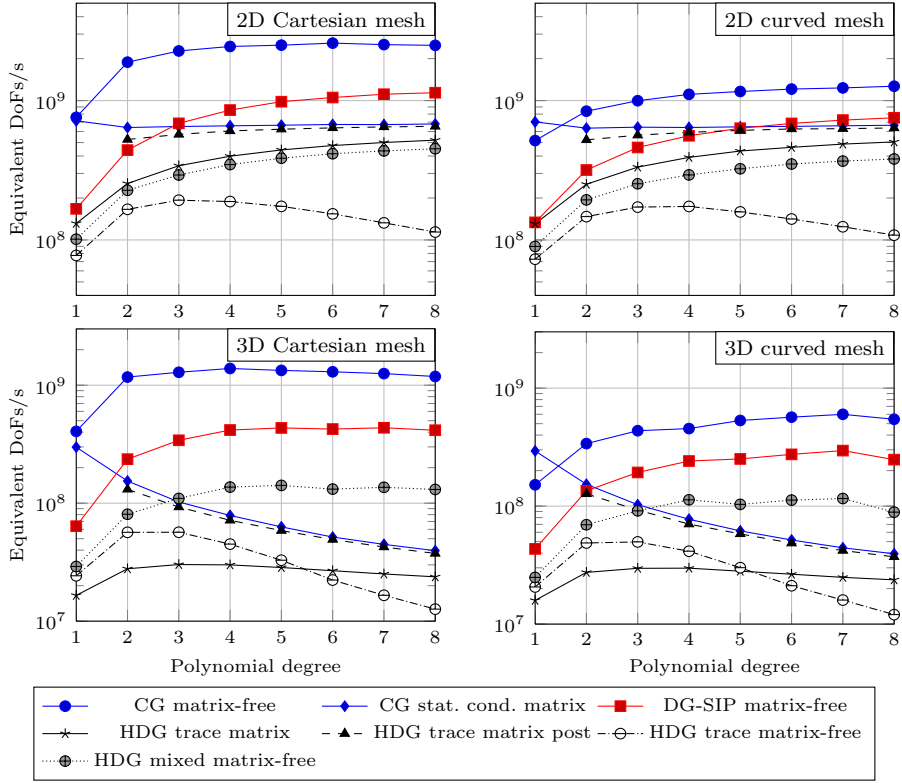

Fig.~\ref{fig:matvec} shows the measurements in two and three space dimensions
on both Cartesian meshes with constant coefficients and curved meshes with
variable coefficients. In all tests, continuous finite elements show
the best performance, apart from linear shape functions.  As a point of
reference, the operator evaluation with our implementation for $\mathcal Q_2$
shape functions on general grids is considerably faster at 340 million degrees
of freedom per second than HPGMG\footnote{\texttt{https://hpgmg.org},
  retrieved on July 20, 2016} at 140 million degrees of freedom per second
($64^3$ mesh), both run on the same hardware. This shows both the high level of
optimization and the benefit of pre-computed Jacobians \cite{kk12}. For
$k\geq 2$, continuous elements provide more than twice the throughput of the
next-best method, the DG-SIP method with matrix-free implementation.  This
goes against the preconception of the DG community \cite{cgl09} which attribute an
implementation advantage to DG methods due to more structure and favorable
vectorization properties. Note that the throughput of the DG-SIP method at 400
million equivalent degrees of freedom per second on the 3D Cartesian mesh
translates to around 500--800 million DG degrees of freedom processed per
second at polynomial degrees two through eight. A remarkable aspect of the
matrix-free schemes is that throughput in terms of degrees of freedom per
second appears approximately constant as the polynomial degree increases,
despite the theoretical $\mathcal O(k)$ complexity per degree of freedom. One
reason for this behavior is that face integrals at cost $\mathcal O(1)$
are dominating at lower degrees $k\leq 4$ for DG-SIP and the
HDG mixed form. Secondly, the even-odd decomposition \cite{Kopriva09} limits
the cost increase at higher degrees. Finally, better arithmetic utilization is
achieved for higher polynomial degrees, as memory access per DoF scales as
$\mathcal O(1)$ in the polynomial degree.

The results also show that HDG with post-processing (black dashed line)
delivers approximately the same performance as matrix-free DG-SIP in two space
dimensions. Without post-processing, HDG falls considerably behind the
matrix-free schemes in 2D. In 3D, all matrix-based schemes are by an order
of magnitude and more
slower than the matrix-free schemes. Our
results show that the post-processed HDG solution at assumed convergence rate $k+2$
performs similarly to the statically condensed CG matrix at rate $k+1$ because both have
the same $k^{d-1}$ unique degrees of freedom per face. In other words, the
effect of superconvergence is offset by the discontinuous trace solution
spaces \eqref{eq:trace_space} and embedded discontinuous Galerkin with traces
based on the skeleton of the continuous finite element method as used in
\cite{npc15} actually appear more favorable in terms of accuracy
efficiency. In the remainder of this study, the statically condensed CG
results can be taken as a proxy for the performance of embedded DG methods.

\subsection{Performance modeling}

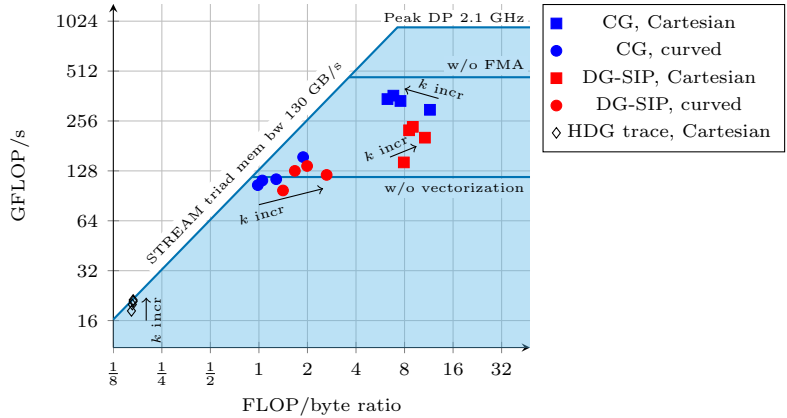
\begin{figure}
\centering
\begin{tikzpicture}
\begin{loglogaxis}
  [scale=0.8,
  grid=both,
  xlabel=FLOP/byte ratio,
  ylabel=GFLOP/s,
  axis lines=left,
  ytick={16,32,64,128,256,512,1024},
  yticklabels={16,32,64,128,256,512,1024},
  xmin=0.125, xmax=48,ymin=11,ymax=1300,
  xtick={0.125,0.25,0.5,1,2,4,8,16,32},
  xticklabels={$\frac 18$,$\frac 14$,$\frac 12$,1,2,4,8,16,32},
  tick label style={font=\scriptsize},
  label style={font=\scriptsize},
  legend style={font=\scriptsize},
  legend pos = outer north east
  ]
  \path[fill=gnuplot@lightblue, fill opacity=0.4] (axis cs:0.125,5) -- (axis cs:0.125,16.25) -- (axis cs:7.23,940) -- (axis cs:48,940) -- (axis cs:48,5);
  \draw[draw=gnuplot@darkblue,thick] (axis cs:0.125,16.25) -- node[above,rotate=46.5,inner sep=2pt,outer sep=1pt,fill=white]{\tiny STREAM triad mem bw 130 GB/s}(axis cs:7.23, 940) -- (axis cs:48, 940) node[anchor=south east,inner sep=1pt,outer sep=1pt,fill=white]{\tiny Peak DP 2.1 GHz};
  \draw[draw=gnuplot@darkblue,thick] (axis cs:3.62, 470) -- (axis cs:48,470) node[anchor=south east,inner sep=1pt,outer sep=1pt,fill=gnuplot@lightblue!40]{\tiny w/o FMA};
  \draw[draw=gnuplot@darkblue,thick] (axis cs:0.904, 117.5) -- (axis cs:48,117.5) node[anchor=north east,inner sep=1pt,outer sep=1pt,fill=gnuplot@lightblue!40]{\tiny w/o vectorization};
  \addplot[only marks,color=blue,mark=square*]
  coordinates{(11.5,299)(6.82,364)(6.27,347)(7.55,338)};
  \addlegendentry{CG, Cartesian};
  \draw [->] (axis cs:13,350) node[above,rotate=-15,inner sep=1pt,outer sep=1pt]{\tiny $k$ incr} -- (axis cs:8,400);
  \addplot[only marks,color=blue,mark=*]
  coordinates{(1.05,112)(0.984,105)(1.28,114)(1.88,155)};
  \addlegendentry{CG, curved};
  \draw [->] (axis cs:1,80) node[below,rotate=15,inner sep=1pt,outer sep=1pt]{\tiny $k$ incr} -- (axis cs:2.5,100);
  \addplot[only marks,color=red,mark=square*]
  coordinates{(7.92,144)(8.54,225)(9.00,236)(10.7,203)};
  \addlegendentry{DG-SIP, Cartesian};
  \draw [->] (axis cs:6.5,155) node[above,rotate=25,inner sep=1pt,outer sep=1pt]{\tiny $k$ incr} -- (axis cs:9.5,180);
  \addplot[only marks,color=red,mark=*]
  coordinates{(1.41,97.5)(1.67,128)(1.99,137)(2.63,121)};
  \addlegendentry{DG-SIP, curved};
  \addplot[only marks,color=black,mark=diamond]
  coordinates{(0.162,18.3)(0.164,20.0)(0.166,20.9)(0.166,21.3)};
  \addlegendentry{HDG trace, Cartesian};
  \draw [->] (axis cs:0.2,16) node[below,rotate=90,inner sep=1pt,outer sep=1pt]{\tiny $k$ incr} -- (axis cs:0.2,22);
\end{loglogaxis}
\end{tikzpicture}
\caption{Roofline model for the evaluation of the 3D Laplacian with different
  variants on a dual-socket Intel Xeon E5-2690 v4 (Broadwell), 28
  cores. Polynomial degrees $k=1, 2, 4, 8$ are considered. Small arrows
  indicate the behavior for increasing polynomial degrees.}
\label{fig:mv_roofline}
\end{figure}

In order to document our unbiased comparison, Fig.~\ref{fig:mv_roofline} puts
the achieved performance into perspective by a roofline performance model. The
performance boundaries are the memory bandwidth limit (diagonal line to the
left) and the arithmetic throughput limit (horizontal line to the upper right)
in terms of the FLOP/byte ratio of the respective kernel
\cite{Patterson09}. The FLOP/byte ratio is measured by an optimistic
assumption to memory access that only counts the global data structures that
need to be streamed at least once per matrix-vector product, assuming perfect
caching of re-usable data such as vector entries accessed by several elements
and no other bottlenecks in the memory hierarchies.
The arithmetic operations per element have been derived
by formulas similar to the ones from \cite{ksmw15} and verified by counting
instructions of a run of the element kernel through the Intel Software
Development Emulator.\footnote{\texttt{https://software.intel.com}, Version
  7.45, AVX2 (Haswell) mode, retrieved on May 19, 2016.}

The results in Fig.~\ref{fig:mv_roofline} show that the sparse matrix-vector
product displayed in the lower left corner is very close to the theoretical
performance maximum of the compressed row storage scheme. The only available
optimization of the matrix-based scheme except for a stencil representation in
the constant-coefficient affine mesh case would be to use the DG matrix
structure where blocks of entries are addressed by one index. In Trilinos, the class
\texttt{Epetra\_VbrMatrix} implements such a scheme rather than the pointwise indirect addressing of the
compressed row storage in \texttt{Epetra\_CrsMatrix}. This would
increase the throughput by up to 50\%,
moving to a FLOP/byte ratio of 0.25. However, inefficiencies and limitations
in the Trilinos implementation prohibit its use in general software such as
a multigrid solver. Even if such an implementation were available, the performance model show today's hardware does not permit matrix-based methods
to match the matrix-free implementations. Turning to the matrix-free implementations,
Fig.~\ref{fig:mv_roofline} reveals a clear difference between the Cartesian
case where the memory transfer is mostly due to the solution vector and some
index data, and the curved mesh case where the Jacobian transformation also
needs to be loaded. The former is computation bound, whereas the latter
resides in the memory bound region, albeit further to the right than sparse
matrix kernels. Due to the considerably more complex kernel structures with
many short loops and different operations, the achieved performance is not as
close to the theoretical performance bounds as for the sparse matrix-vector
kernel. Given that our implementation clearly outperforms the benchmark code
HPGMG, the numbers can be considered extremely good nonetheless.

\subsection{Latency analysis}

\begin{figure}
  \centering
\pgfplotstableread{
        nprocs     femf       dgmf       spmv
        16         0.00931212 0.0433836          NaN
        32         0.00464533 0.0219271          NaN
        64         0.00230463 0.010921    3.0112e-02
        128        0.00117246 0.00669329  1.5112e-02
        256        0.00063825 0.00341304  7.8186e-03
        512        0.00035622 0.00173576  4.1118e-03
        1024       0.00021754 0.00092641  2.2146e-03
        2048       0.00014925 0.00051910  1.2636e-03
        4096       0.00012050 0.00030556  7.5517e-04
        8192       0.00011000 0.00019982  5.6549e-04
        16384      0.00010985 0.00016806  5.8049e-04
        32768      0.00016206 0.00021245  4.7015e-04
        65536      9.8138e-05 0.00012694  4.6373e-04
}\mytable
  \begin{tikzpicture}
    \begin{loglogaxis}[
      width=0.55\textwidth,
      height=0.42\textwidth,
      xlabel={Number of cores},
      ylabel={Time matrix-vector product [s]},
      x label style={at={(0.5,0.05)}},
      tick label style={font=\scriptsize},
      label style={font=\scriptsize},
      legend style={font=\scriptsize},
      xtick={16,64,256,1024,4096,16384,65536},
      xticklabels={16,64,256,1024,4096,16k,64k},
      ymin=5e-5,
      ymax=1e-1,
      grid
      ]
      \addplot table[x={nprocs}, y={femf}] {\mytable};
      \addlegendentry{CG mat-free};
      \addplot table[x={nprocs}, y={dgmf}] {\mytable};
      \addlegendentry{DG-SIP mat-free};
      \pgfplotsset{cycle list shift=1}
      \addplot table[x={nprocs}, y={spmv}] {\mytable};
      \addlegendentry{HDG matrix};
    \end{loglogaxis}
\end{tikzpicture}
\caption{Latency study of matrix-vector products for Laplacian on a $80^3$
  mesh with $\mathcal Q_2$ elements, involving $4.1$ million DoFs.}
\label{fig:matvec_latency}
\end{figure}
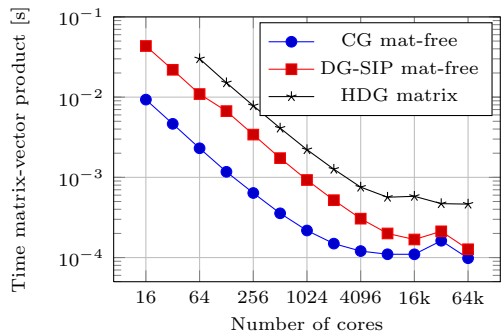

A second ingredient to practical solver performance is the latency of matrix-vector
products. This is relevant for multigrid schemes where a series of coarser
representations appear and need to be processed quickly. This section reports
results from the SuperMUC Phase 1 system ($2\times 8$ core Intel Xeon E5-2680
Sandy Bridge CPU running at 2.7~GHz, Infiniband FDR10
interconnect). Fig.~\ref{fig:matvec_latency} shows that the matrix-free
variants scale down to $10^{-4}$ seconds where network latency becomes
dominant, a similar result as was reported for HPGMG \cite{hpgmg16}. This number
needs to be compared to a single point-to-point latency of around $10^{-6}$
seconds.

The sparse matrix-vector product of HDG already saturates at around
$5\cdot 10^{-4}$ seconds. Investigation of this issue revealed sub-optimal MPI
commands in the data exchange routines of the Epetra sparse matrix
\cite{trilinos}, involving a global barrier operation in addition to
point-to-point communication. However, the latency could not be reduced to
less than $2\cdot 10^{-4}$ seconds even when changing the Trilinos source code.
We conclude that given the right implementation, no
advantage of the discontinuous data structures with less connectivity to
neighbors appears on modern high-performance implementations of matrix-vector
products, as opposed to direct solvers \cite{harp13}, showing that the results
from \cite{Yakovlev16} are not general. When going to considerably higher
order elements than the $\mathcal Q_2$ basis reported in Fig.~\ref{fig:matvec_latency},
the time to process a single element, or rather, a batch of four to eight
elements due to vectorization over several elements in our implementation
\cite{kk12}, overlays the communication latency. For DG-SIP, this transition
occurs at $k=6$ and at around $k=8$ for CG.

\section{Performance comparison of multigrid solvers}\label{sec:comparison}

In this section, we analyze modern multigrid solvers and record the solution
accuracy as a function of computing time. All solvers
use a multigrid V-cycle as a preconditioner for a conjugate gradient
iteration, which increases solver robustness \cite{tos01,Sundar15}. The
iteration is stopped once the residual norm
goes below $10^{-9}$ times the right hand side norm. For coarse and low-order
discretizations, this tolerance could be relaxed as e.g.~done by the full
multigrid cycle in HPGMG \cite{hpgmg16}, but we refrain from this optimization
for ease of comparison. The comparisons focus on the three-dimensional case
where large-scale iterative solvers are essential. The trend for 2D is similar
to 3D, but the advantage of the matrix-free schemes is less. In 2D, they
provide two to eight times better efficiency, solving for 10 and 20 million
unknowns per seconds on 28 cores for the continuous and discontinuous
matrix-free solvers. The HDG solver reaches between 1 and 3 million degrees of
freedom per second.

\subsection{Multigrid solvers for matrix-free methods}\label{sec:mg_cg}

Due to their fast matrix-vector products, a polynomial Chebyshev accelerated
pointwise Jacobi smoother in a geometric multigrid cycle is the natural choice
\cite{ABHT03,hpgmg16,kk12} for the CG and DG-SIP realizations. Besides
matrix-vector products, the Chebyshev smoother only needs access to the matrix
diagonal that is pre-computed and stored before solving. In addition, an
estimate of the largest eigenvalue $\tilde{\lambda}_{\max{}}$ of the
Jacobi-preconditioned matrix is used to make the Chebyshev iteration address
modes with eigenvalue in the interval
$[0.06\tilde{\lambda}_{\max{}},1.2\tilde{\lambda}_{\max{}}]$. The eigenvalue
estimation is done by 15 iterations with the conjugate gradient method. In the
numbers reported below, the setup cost of the multigrid ingredients is
ignored, just as we ignore the cost for assembling the HDG trace matrix. We
note that the setup of the matrix-free variants is proportional to two to four
V-cycles, considerably less than the matrix creation and assembly for HDG. In
the context of nonlinear systems where the system matrix changes rapidly, the
advantage of matrix-free schemes will thus be even larger than what is
reported here, a property exploited in \cite{Kronbichler16}.

For pre- and post-smoothing, a polynomial degree of five in the Chebyshev
method is used, involving five matrix-vector products. The level transfer is
based on the usual geometric embedding operations and also implemented by
tensorial techniques.  For the solver on the coarse grid, the Chebyshev
iteration is selected, now with parameters such the a-priori error estimate
for the Chebyshev iteration \cite{Varga09} ensures an error below
$10^{-3}$. The implementation
uses the multigrid facilities of the \texttt{deal.II} finite element library
\cite{jk11,Kanschat04}, including adaptively refined meshes with
hanging nodes in a massively parallel context based on a forest-of-tree data
layout and Morton cell ordering \cite{bbhk11,p4est}.

\subsection{Multigrid for HDG}

In the context of the HDG trace system, off-the-shelf AMG solvers such as
Trilinos ML \cite{ml-guide} work suboptimally or even fail because of a
pronounced non-diagonally dominant character of the matrix together with wide
stencils due to the high-order basis. To overcome these limitations, this work
adopts a variation of the method proposed in \cite{cdgt13}, where the
combination of a high-order HDG trace space with a continuous finite element
discretization involving linear basis functions \emph{on the same mesh} was
proposed. This concept is closely related to $p$-multigrid methods
\cite{lbl06} where the structure of a high-order basis is used by first going to a low-order basis with fewer unknowns
rather than going to coarser meshes as in $h$-multigrid
approaches. The transfer between the HDG trace space and linear finite elements is
realized by the embedding operator that maps linear shape functions onto the
trace polynomials as well as its transpose. As opposed to the work \cite{cdgt13}
that constructs a genuine discretization on the linear finite element space, we select
a Galerkin coarse grid operator \cite{tos01}. As shown below, the iteration
count is only around 15--20, much better than 55--75 reported in \cite{cdgt13}
for similar tolerances that are probably to systematic gaps between different
discretizations. The hierarchy is then continued by algebraic multigrid. Since
the connectivity of the matrix is the same as for linear finite elements and
the matrix is (almost) diagonally dominant, optimal or close-to-optimal
performance of the AMG inside the $p$-AMG scheme can be expected.

Due to a strong non-diagonally dominant matrix structure, optimal multigrid
performance in HDG cannot be obtained with point-relaxation smoothers. Instead,
block-relaxation scheme with blocks combining all
degrees of freedom on a face or incomplete factorizations are
necessary. (Iteration numbers grow approximately as $h^{-2/3}$ with point
Gauss--Seidel smoothing.) Due to its robustness, ILU(0) is selected as a
smoother both for the HDG trace matrix as well as in the AMG levels. As soon
as the level matrix size goes below 2000, a direct coarse solver is
invoked. For the HDG stabilization parameter, we select
$\tau = 5 \max\{\kappa(\vec x), \vec x \in \text{face}\}$ as a balance between
solver efficiency and accuracy throughout this study, see also \cite{ksc12}.

For the HDG trace matrix, we consider the representation by a sparse matrix
because it easily combines with the ILU(0) for the smoother, even though
somewhat higher performance would be available in 3D for $1\leq k\leq 5$ with
a matrix-free implementation according to Fig.~\ref{fig:matvec}. Operator
evaluation with the mixed form of HDG seems promising due to the considerably
faster matrix-vector product reported in Fig.~\ref{fig:matvec}, but the saddle
point form is more challenging to handle, requiring strong ingredients such as
overlapping Schwarz smoothers or block factorizations \cite{bgl05}.

For comparison, we also consider an iterative solver based on the statically
condensed matrix for continuous elements. The standard ML-AMG V-cycle with one
sweep of ILU(0) for pre- and post-smoothing on all levels is chosen. As seen
from Table~\ref{tab:solver_times_3d} below, no optimal iteration numbers are
obtained in this case for higher polynomial degrees. However, it serves as a
point of reference for matrix-based approaches with black-box
preconditioning. Alternatively, a similar $p$-multigrid scheme with similar
iteration counts as for the HDG trace system could also be considered.

\subsection{Three-dimensional example with smooth solution}

\begin{figure}
  \centering
  \includegraphics[width=0.45\textwidth]{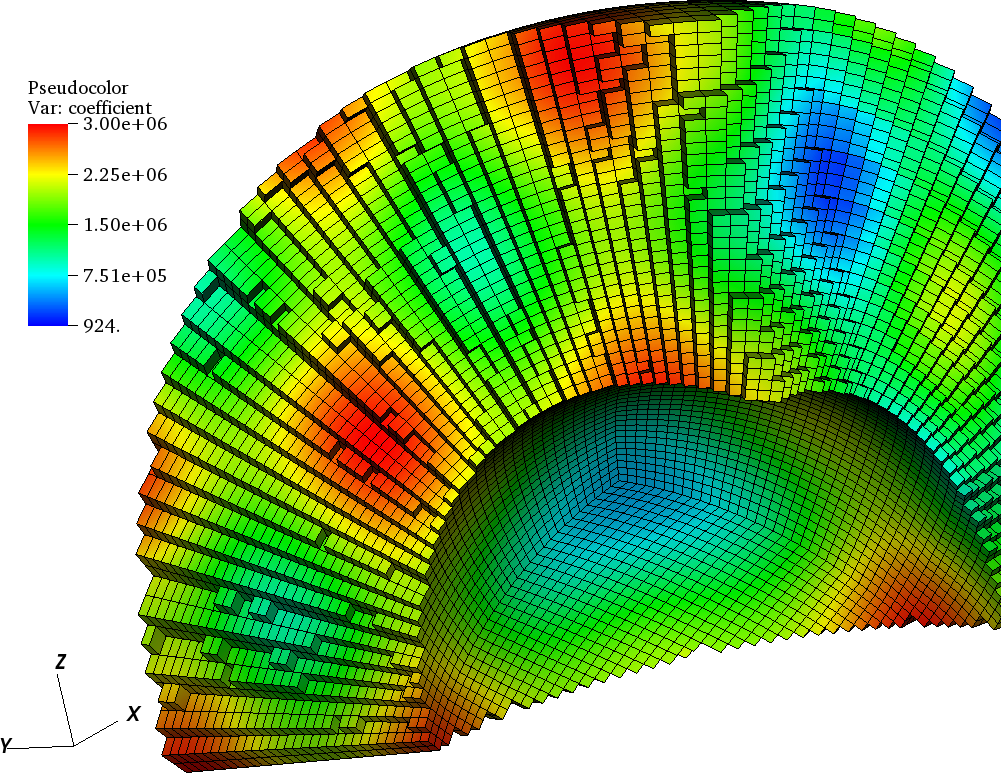}
  \caption{Visualization of parts of the three-dimensional shell geometry with
    the variable coefficient \eqref{eq:diffusivity_shell} on a mesh consisting
    of 393k elements in total.}
  \label{fig:mesh_shell_3d}
\end{figure}

We consider the Poisson equation with analytic solution
\begin{equation}\label{eq:pois_analytic}
u(\vec x) = \left(\frac{1}{\alpha \sqrt{2\pi}}\right)^3\sum_{j=1}^3 \mathrm{exp}\left(- \|\vec x - \vec x_j\|^2/\alpha ^2\right),
\end{equation}
given as a sum of three Gaussians centered at the positions
$\vec x_j \in \{(-0.5, 0.5, 0.25)^\mathrm T,\allowbreak (-0.6,-0.5,-0.125)^\mathrm{T},\allowbreak (0.5,-0.5,0.5)^\mathrm T \}$
and of width $\alpha=\frac 15$.  The equation is solved on two domains,
\begin{itemize}
\item the unit cube, $\Omega = (-1,1)^3$, with the surfaces at $x_e=-1$ subject to
  Neumann boundary conditions and the surfaces at $x_e=+1$ subject to Dirichlet
  boundary conditions, $e=1,2,3$, using constant diffusivity $\kappa = 1$, and
\item a full spherical shell in 3D with inner radius $0.5$ and outer radius $1.0$,
  using a polynomial approximation of degree 5 along a spherical
  manifold for all elements, and a strongly varying diffusivity
  \begin{equation}\label{eq:diffusivity_shell}
    \kappa(\vec x) = \kappa (x_1, \ldots, x_d) =  1 + 10^6 \prod_{e=1}^3 \cos(2\pi x_e + 0.1 e),
  \end{equation}
  inspired by the variable coefficient case in \cite{Sundar15} but using a
  shift $0.1e$ in order to eliminate any potential spatial symmetries in the
  coefficients. Dirichlet conditions are set on all boundaries. A
  visualization of approximately one eighth of a sample 3D mesh along with the
  coefficient is given in Fig.~\ref{fig:mesh_shell_3d}. Geometric multigrid
  methods start with an initial mesh consisting of six elements.
\end{itemize}

In both cases, the boundary conditions $g_\text{D}$ and $g_\text{N}$ as well
as the value of the forcing $f$ in \eqref{eq:poisson} are set such that the
analytic solution \eqref{eq:pois_analytic} is obtained.

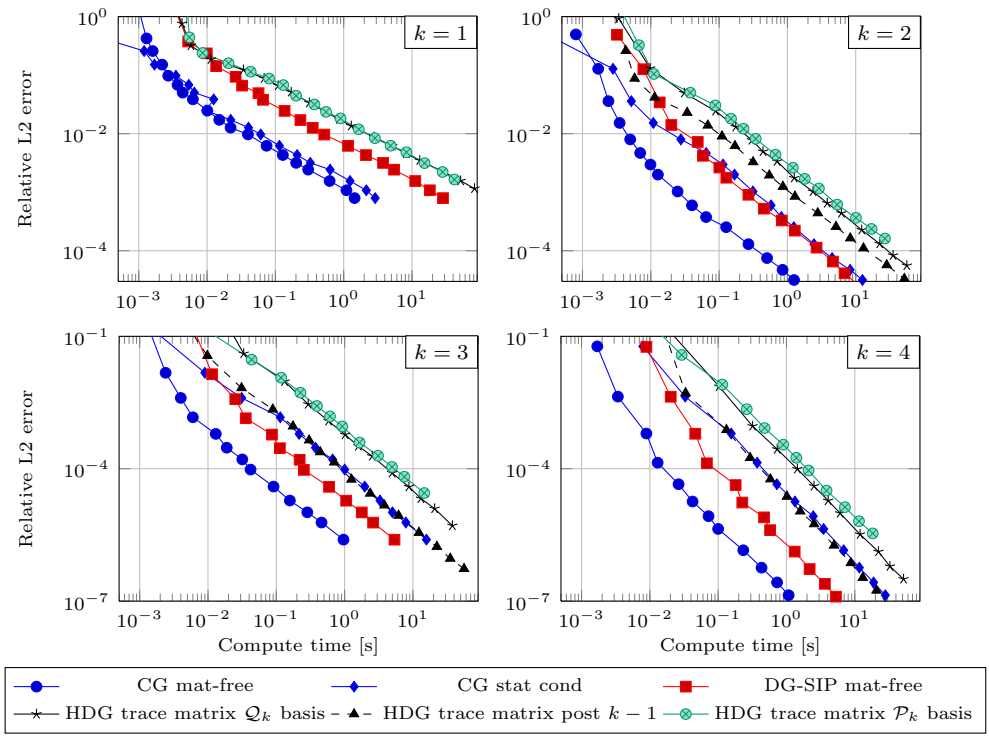
\begin{figure}
  \centering
\pgfplotstableread{
 cells    Gtime  Gits    hdgL2 hdgL2post  qTime  qIt    FEL2  trTime trIt  dgTim dgIt dgL2     hdgpTi  hdgpIt hdgpEr
      64  0.0022  12 3.874e+00 3.519e+00  0.0008 4 3.517e+00 0.0001   1  0.0024  14  3.468e+00  0.0019  12 1.649e+00
     512  0.0043  13 7.612e-01 2.601e-01  0.0013 4 4.241e-01 0.0003   1  0.0053  15  3.759e-01  0.0051  14 1.566e+00
    1728  0.0058  13 3.189e-01 8.775e-02  0.0016 4 2.580e-01 0.0012   9  0.0099  14  2.354e-01  0.0055  14 4.422e-01
    4096  0.0114  13 1.850e-01 4.167e-02  0.0022 4 1.511e-01 0.0017  10  0.0136  14  1.422e-01  0.0086  14 2.385e-01
    8000  0.0341  13 1.232e-01 2.286e-02  0.0027 4 9.768e-02 0.0035   9  0.0262  14  9.373e-02  0.0206  14 1.598e-01
   13824  0.0684  13 8.799e-02 1.385e-02  0.0037 4 6.821e-02 0.0054   9  0.0326  14  6.622e-02  0.0435  14 1.145e-01
   21952  0.1102  13 6.598e-02 9.015e-03  0.0044 4 5.029e-02 0.0065   9  0.0566  14  4.918e-02  0.0809  13 8.740e-02
   32768  0.1706  13 5.131e-02 6.190e-03  0.0062 4 3.859e-02 0.0125   9  0.0660  14  3.793e-02  0.1305  14 6.758e-02
   64000  0.3146  14 3.357e-02 3.281e-03  0.0101 4 2.477e-02 0.0096  11  0.1364  14  2.449e-02  0.2009  14 4.453e-02
  110592  0.5432  14 2.366e-02 1.944e-03  0.0151 4 1.723e-02 0.0222  12  0.2326  14  1.709e-02  0.3729  15 3.155e-02
  175616  0.8607  14 1.757e-02 1.245e-03  0.0221 4 1.267e-02 0.0404  11  0.3521  14  1.259e-02  0.5522  14 2.353e-02
  262144  1.2898  14 1.356e-02 8.446e-04  0.0397 5 9.704e-03 0.0610  10  0.5200  14  9.661e-03  0.8979  14 1.823e-02
  512000  2.7820  15 8.780e-03 4.403e-04  0.0743 5 6.215e-03 0.1147  10  1.1578  14  6.197e-03  1.6621  14 1.189e-02
  884736  5.1425  16 6.144e-03 2.579e-04  0.1286 5 4.318e-03 0.2083  10  2.1217  14  4.309e-03  2.8621  14 8.374e-03
 1404928  8.3169  16 4.539e-03 1.638e-04  0.2015 5 3.173e-03 0.3245  10  3.7677  14  3.168e-03  4.8667  15 6.227e-03
 2097152 12.9887  17 3.490e-03 1.105e-04  0.3024 5 2.430e-03 0.6326  13  5.4621  14  2.427e-03  8.2908  17 4.818e-03
 4096000 28.4178  19 2.247e-03 5.710e-05  0.6225 5 1.555e-03 1.2385  13 11.1565  14  1.554e-03 15.2781  16 3.140e-03
 7077888 51.5241  20 1.566e-03 3.325e-05  1.0970 5 1.080e-03 2.1432  12 18.2271  14  1.080e-03 28.0464  17 2.216e-03
11239424 82.1699  20 1.154e-03 2.103e-05  1.4438 4 7.937e-04 2.8896  11 28.4151  14  7.933e-04 41.8153  16 1.651e-03
}\errorsdA
\pgfplotstableread{
 cells    Gtime  Gits    hdgL2 hdgL2post  qTime  qIt     FEL2  trTime trIt  dgTim dgIt dgL2     hdgpTi hdgpIt hdgpEr
     64  0.0034   13 9.250e-01 4.227e-01 0.0008  4   4.958e-01 0.0003   1   0.0032 12 4.893e-01  0.0025 15 3.105e+00
    512  0.0098   14 1.307e-01 3.628e-02 0.0017  4   1.283e-01 0.0028  14   0.0078 12 1.264e-01  0.0067 15 3.249e-01
   1728  0.0310   14 4.950e-02 6.753e-03 0.0024  4   3.595e-02 0.0052  13   0.0135 11 3.424e-02  0.0111 14 1.047e-01
   4096  0.0888   14 2.434e-02 2.216e-03 0.0035  4   1.526e-02 0.0109  14   0.0200 11 1.409e-02  0.0379 15 5.013e-02
   8000  0.1757   14 1.305e-02 9.151e-04 0.0050  4   7.982e-03 0.0276  15   0.0487 12 7.235e-03  0.0881 14 3.043e-02
  13824  0.3023   14 7.741e-03 4.420e-04 0.0070  4   4.678e-03 0.0648  15   0.0580 12 4.193e-03  0.1496 14 1.805e-02
  21952  0.4394   13 4.960e-03 2.386e-04 0.0099  4   2.968e-03 0.1142  16   0.1009 11 2.641e-03  0.2277 14 1.215e-02
  32768  0.7059   14 3.365e-03 1.398e-04 0.0127  4   1.998e-03 0.1685  16   0.1277 12 1.769e-03  0.3436 14 8.189e-03
  64000  1.2591   13 1.753e-03 5.715e-05 0.0247  4   1.029e-03 0.3132  18   0.2683 12 9.055e-04  0.6681 14 4.386e-03
 110592  2.3699   14 1.026e-03 2.752e-05 0.0401  4   5.973e-04 0.5766  19   0.4535 12 5.238e-04  1.2036 15 2.622e-03
 175616  3.8289   14 6.514e-04 1.483e-05 0.0645  4   3.768e-04 0.8158  17   0.8270 12 3.298e-04  1.7998 14 1.693e-03
 262144  6.1754   15 4.389e-04 8.686e-06 0.1269  5   2.528e-04 1.2437  17   1.2825 12 2.209e-04  2.8895 15 1.159e-03
 512000 12.2448   15 2.266e-04 3.552e-06 0.2691  5   1.296e-04 2.4510  17   2.6744 12 1.131e-04  5.3556 14 6.143e-04
 884736 22.4358   16 1.318e-04 1.711e-06 0.5032  5   7.506e-05 4.5180  18   4.6393 12 6.543e-05 10.0777 15 3.659e-04
1404928 35.2064   16 8.334e-05 9.228e-07 0.8499  5   4.729e-05 8.2848  21   6.8660 11 4.120e-05 17.0513 16 2.362e-04
2097152 56.0182   17 5.599e-05 5.406e-07 1.2501  5   3.169e-05 12.5648 21  10.2283 11 2.760e-05 26.7925 17 1.619e-04
}\errorsdB
\pgfplotstableread{
 cells    Gtime  Gits    hdgL2 hdgL2post  qTime  qIt     FEL2  trTime trIt  dgTim dgIt dgL2    hdgpTi hdgpIt hdgpEr
    64  0.0161   14  2.931e-01 1.360e-01  0.0012  4 2.511e-01  0.0010  1  0.0054  11 2.481e-01  0.0043 15  8.002e-01
   512  0.0333   16  4.061e-02 5.209e-03  0.0024  4 1.491e-02  0.0090 17  0.0115  11 1.385e-02  0.0127 16  1.027e-01
  1728  0.1321   16  9.495e-03 7.597e-04  0.0040  4 4.023e-03  0.0311 18  0.0251  10 3.770e-03  0.0435 15  2.964e-02
  4096  0.2934   16  2.894e-03 1.767e-04  0.0060  4 1.465e-03  0.1143 19  0.0357  10 1.394e-03  0.1178 15  1.153e-02
  8000  0.5814   16  1.203e-03 5.789e-05  0.0130  5 6.155e-04  0.2158 19  0.0863  10 5.937e-04  0.2256 15  5.284e-03
 13824  1.0168   16  5.926e-04 2.339e-05  0.0188  5 2.997e-04  0.3745 19  0.1127  10 2.917e-04  0.3937 15  2.638e-03
 21952  1.5900   16  3.247e-04 1.085e-05  0.0318  5 1.627e-04  0.6673 23  0.2197  10 1.594e-04  0.6109 15  1.535e-03
 32768  2.5370   17  1.924e-04 5.574e-06  0.0420  5 9.572e-05  1.0039 23  0.2529  10 9.419e-05  0.9232 15  9.087e-04
 64000  4.9630   17  7.994e-05 1.829e-06  0.0911  5 3.938e-05  1.9710 23  0.5919  10 3.896e-05  1.6534 14  3.927e-04
110592  8.7499   17  3.891e-05 7.354e-07  0.1578  5 1.904e-05  3.2953 22  1.0542  10 1.889e-05  3.0849 15  1.969e-04
175616 13.0962   16  2.114e-05 3.402e-07  0.2845  5 1.029e-05  5.0106 21  1.7842  10 1.023e-05  4.9391 15  1.096e-04
262144 20.7953   17  1.245e-05 1.745e-07  0.4598  5 6.038e-06  7.8852 22  2.6152  10 6.012e-06  7.4939 15  6.590e-05
512000 37.8215   16  5.132e-06 5.716e-08  0.9624  5 2.476e-06 15.7584 22  5.3714  10 2.469e-06 14.6613 15  2.812e-05
}\errorsdC
\pgfplotstableread{
 cells    Gtime  Gits   hdgL2 hdgL2post  qTime  qIt    FEL2  trTime trIt  dgTim dgIt dgL2    hdgpTi hdgpIt hdgpEr
    64  0.0201   15 1.235e-01 3.437e-02  0.0017 4 5.924e-02  0.0079  14 0.0089 11 5.760e-02  0.0079  15 2.809e-01
   512  0.1034   17 7.295e-03 7.857e-04  0.0034 4 4.305e-03  0.0329  20 0.0202 11 4.271e-03  0.0292  17 3.787e-02
  1728  0.3192   17 9.237e-04 6.783e-05  0.0089 5 6.270e-04  0.1553  21 0.0462 10 6.202e-04  0.1147  16 8.081e-03
  4096  0.7309   17 2.767e-04 1.340e-05  0.0130 5 1.370e-04  0.3727  22 0.0684 10 1.331e-04  0.2617  16 2.229e-03
  8000  1.4454   17 9.859e-05 3.661e-06  0.0263 5 4.457e-05  0.7237  22 0.1803 10 4.267e-05  0.4787  15 8.350e-04
 13824  2.5294   17 4.043e-05 1.231e-06  0.0422 5 1.807e-05  1.3327  23 0.2257 10 1.714e-05  0.8942  16 3.507e-04
 21952  4.0742   17 1.895e-05 4.890e-07  0.0722 5 8.405e-06  2.4664  27 0.4697 10 7.919e-06  1.3955  16 1.768e-04
 32768  6.0984   17 9.814e-06 2.196e-07  0.0999 5 4.326e-06  3.5002  26 0.5759 10 4.058e-06  2.1026  16 9.158e-05
 64000 11.9522   17 3.257e-06 5.758e-08  0.2331 5 1.423e-06  6.8847  26 1.3182 10 1.328e-06  3.8581  15 3.198e-05
110592 22.2282   18 1.319e-06 1.928e-08  0.4315 5 5.733e-07 11.6710  25 2.1823 10 5.330e-07  7.1735  16 1.347e-05
175616 32.6582   17 6.137e-07 7.640e-09  0.7286 5 2.656e-07 18.8503  25 3.6454 10 2.464e-07 11.5196  16 6.470e-06
262144 51.7283   18 3.160e-07 3.427e-09  1.0753 5 1.363e-07 27.9945  25 5.2967 10 1.263e-07 18.2992  17 3.424e-06
}\errorsdD
\pgfplotstableread{
 cells    Gtime  Gits   hdgL2 hdgL2post  qTime  qIt    FEL2  trTime trIt  dgTim dgIt dgL2    hdgpTi hdgpIt hdgpEr
    64  0.0445   16 4.983e-02 9.292e-03 0.0033  5 2.366e-02  0.0218  17 0.0189  13 2.344e-02  0.0161 17 1.437e-01
   512  0.2189   18 1.368e-03 1.196e-04 0.0064  5 5.964e-04  0.1170  22 0.0379  12 5.734e-04  0.0727 18 1.453e-02
  1728  0.6759   18 1.707e-04 8.627e-06 0.0142  5 4.680e-05  0.3988  23 0.1053  12 4.431e-05  0.2327 17 1.909e-03
  4096  1.6534   19 2.974e-05 1.123e-06 0.0240  5 1.067e-05  0.9183  23 0.1525  12 1.022e-05  0.5224 17 3.990e-04
  8000  3.3229   19 7.364e-06 2.230e-07 0.0523  5 3.013e-06  1.8749  24 0.3774  11 2.911e-06  1.0271 17 1.233e-04
 13824  5.8154   19 2.508e-06 6.244e-08 0.0822  5 1.018e-06  3.8348  28 0.5079  11 9.922e-07  1.8095 17 4.245e-05
 21952  9.3104   19 1.007e-06 2.127e-08 0.1582  5 4.060e-07  5.8964  27 1.1683  11 3.979e-07  2.8634 17 1.876e-05
 32768 13.9873   19 4.558e-07 8.361e-09 0.2277  5 1.828e-07  9.1800  28 1.3863  12 1.799e-07  4.3085 17 8.391e-06
 64000 27.1941   19 1.208e-07 1.756e-09 0.5442  5 4.812e-08 16.9694  26 2.7474  11 4.761e-08  7.8820 16 2.364e-06
110592 49.1285   20 4.076e-08 4.929e-10 0.9218  5 1.615e-08 29.4182  26 4.4418  11 1.603e-08 14.7362 17 8.364e-07
}\errorsdE
  \begin{tikzpicture}
    \begin{loglogaxis}[
      title style={at={(1,0.944)},anchor=north east,draw=black,fill=white,font=\footnotesize},
      title={$k=1$},
      width=0.52\textwidth,
      height=0.42\textwidth,
      ylabel={Relative L2 error},
      x label style={at={(0.5,0.05)}},
      tick label style={font=\scriptsize},
      label style={font=\scriptsize},
      legend style={font=\scriptsize},
      xmin=5e-4, xmax=9e1,
      ymin=3e-5, ymax=1e-0,
      grid
      ]
      \addplot table[x={qTime}, y={FEL2}] {\errorsdA};
      \pgfplotsset{cycle list shift=3}
      \addplot table[x={trTime}, y={FEL2}] {\errorsdA};
      \pgfplotsset{cycle list shift=-1}
      \addplot table[x={dgTim}, y={dgL2}] {\errorsdA};
      \pgfplotsset{cycle list shift=0}
      \addplot table[x={Gtime}, y={hdgL2}] {\errorsdA};
      \addplot[gnuplot@green,mark=otimes*,mark options={fill=gnuplot@green!40}] table[x={hdgpTi}, y={hdgpEr}] {\errorsdA};
    \end{loglogaxis}
  \end{tikzpicture}
  \
  \begin{tikzpicture}
    \begin{loglogaxis}[
      title style={at={(1,0.944)},anchor=north east,draw=black,fill=white,font=\footnotesize},
      title={$k=2$},
      width=0.52\textwidth,
      height=0.42\textwidth,
      x label style={at={(0.5,0.05)}},
      tick label style={font=\scriptsize},
      label style={font=\scriptsize},
      legend style={font=\scriptsize},
      legend to name = legend3dC,
      legend columns = 3,
      xmin=5e-4, xmax=9e1,
      ymin=3e-5, ymax=1e0,
      grid
      ]
      \addplot table[x={qTime}, y={FEL2}] {\errorsdB};
      \addlegendentry{CG mat-free};
      \pgfplotsset{cycle list shift=3}
      \addplot table[x={trTime}, y={FEL2}] {\errorsdB};
      \addlegendentry{CG stat cond};
      \pgfplotsset{cycle list shift=-1}
      \addplot table[x={dgTim}, y={dgL2}] {\errorsdB};
      \addlegendentry{DG-SIP mat-free};
      \pgfplotsset{cycle list shift=0}
      \addplot table[x={Gtime}, y={hdgL2}] {\errorsdB};
      \addlegendentry{HDG trace matrix $\mathcal Q_k$ basis};
      \addplot[mark=triangle*,dashed,black,every mark/.append style={solid}] table[x={Gtime}, y={hdgL2post}] {\errorsdA};
      \addlegendentry{HDG trace matrix post $k-1$};
      \addplot[gnuplot@green,mark=otimes*,mark options={fill=gnuplot@green!40}] table[x={hdgpTi}, y={hdgpEr}] {\errorsdB};
      \addlegendentry{HDG trace matrix $\mathcal P_k$ basis};
    \end{loglogaxis}
  \end{tikzpicture}
  \\
  \begin{tikzpicture}
    \begin{loglogaxis}[
      title style={at={(1,0.944)},anchor=north east,draw=black,fill=white,font=\footnotesize},
      title={$k=3$},
      width=0.52\textwidth,
      height=0.42\textwidth,
      xlabel={Compute time [s]},
      ylabel={Relative L2 error},
      x label style={at={(0.5,0.05)}},
      tick label style={font=\scriptsize},
      label style={font=\scriptsize},
      legend style={font=\scriptsize},
      xmin=5e-4, xmax=9e1,
      ymin=1e-7, ymax=1e-1,
      grid
      ]
      \addplot table[x={qTime}, y={FEL2}] {\errorsdC};
      \pgfplotsset{cycle list shift=3}
      \addplot table[x={trTime}, y={FEL2}] {\errorsdC};
      \pgfplotsset{cycle list shift=-1}
      \addplot table[x={dgTim}, y={dgL2}] {\errorsdC};
      \pgfplotsset{cycle list shift=0}
      \addplot table[x={Gtime}, y={hdgL2}] {\errorsdC};
      \addplot[mark=triangle*,dashed,black,every mark/.append style={solid}] table[x={Gtime}, y={hdgL2post}] {\errorsdB};
      \addplot[gnuplot@green,mark=otimes*,mark options={fill=gnuplot@green!40}] table[x={hdgpTi}, y={hdgpEr}] {\errorsdC};
    \end{loglogaxis}
  \end{tikzpicture}
  \
  \begin{tikzpicture}
    \begin{loglogaxis}[
      title style={at={(1,0.944)},anchor=north east,draw=black,fill=white,font=\footnotesize},
      title={$k=4$},
      width=0.52\textwidth,
      height=0.42\textwidth,
      xlabel={Compute time [s]},
      x label style={at={(0.5,0.05)}},
      tick label style={font=\scriptsize},
      label style={font=\scriptsize},
      legend style={font=\scriptsize},
      xmin=5e-4, xmax=9e1,
      ymin=1e-7, ymax=1e-1,
      grid
      ]
      \addplot table[x={qTime}, y={FEL2}] {\errorsdD};
      \pgfplotsset{cycle list shift=3}
      \addplot table[x={trTime}, y={FEL2}] {\errorsdD};
      \pgfplotsset{cycle list shift=-1}
      \addplot table[x={dgTim}, y={dgL2}] {\errorsdD};
      \pgfplotsset{cycle list shift=0}
      \addplot table[x={Gtime}, y={hdgL2}] {\errorsdD};
      \addplot[mark=triangle*,dashed,black,every mark/.append style={solid}] table[x={Gtime}, y={hdgL2post}] {\errorsdC};
      \addplot[gnuplot@green,mark=otimes*,mark options={fill=gnuplot@green!40}] table[x={hdgpTi}, y={hdgpEr}] {\errorsdD};
    \end{loglogaxis}
  \end{tikzpicture}
\\
\ref{legend3dC}
\caption{Accuracy over solver time for 3D Laplacian on the unit cube.}
\label{fig:accuracy_3d_cube}
\end{figure}

\begin{figure}
  \centering
\pgfplotstableread{
 cells   Gtime  Gits     hdgL2 hdgL2post  qTime  qIt    FEL2  trTime trIt  dgTim dgIt dgL2
     12   0.0010   8 1.640e+00 1.529e+00 0.0003  0 1.516e+00 0.0000  0     0.0007  3 1.327e+00
     96   0.0016  11 1.363e+00 1.004e+00 0.0006  4 1.021e+00 0.0001  1     0.0030 14 9.088e-01
    768   0.0089  19 3.859e-01 1.643e-01 0.0013  4 3.312e-01 0.0003  1     0.0067 15 2.687e-01
   6144   0.0417  27 1.154e-01 6.337e-02 0.0035  5 9.156e-02 0.0021 11     0.0201 14 7.603e-02
  49152   0.8741  51 3.414e-02 4.266e-03 0.0137  5 2.350e-02 0.0074 12     0.1241 13 1.983e-02
 393216  30.4898 217 9.527e-03 6.611e-04 0.1225  5 5.924e-03 0.1179 13     1.1278 13 5.028e-03
3145728 978.0592 839 2.568e-03 1.068e-04 1.0187  5 1.484e-03 0.9431 14    10.0099 13 1.262e-03
}\errorsdA
\pgfplotstableread{
 cells    Gtime  Gits    hdgL2 hdgL2post  qTime  qIt    FEL2  trTime trIt  dgTim dgIt dgL2
     12  0.0040  8  2.566e+00 1.782e+00  0.0004  3 1.754e+00  0.0002   1  0.0020  3 1.769e+00
     96  0.0074 12  5.115e-01 2.622e-01  0.0009  4 3.414e-01  0.0003   1  0.0048 12 2.897e-01
    768  0.0301 16  8.874e-02 2.736e-02  0.0019  4 5.363e-02  0.0030  14  0.0104 12 4.933e-02
   6144  0.0960 19  1.575e-02 1.836e-02  0.0052  4 7.105e-03  0.0237  16  0.0365 12 6.468e-03
  49152  0.7655 24  2.333e-03 1.376e-04  0.0349  4 9.292e-04  0.2653  20  0.3216 12 8.349e-04
 393216  6.1062 32  3.321e-04 2.092e-05  0.3892  4 1.174e-04  2.3345  21  2.9419 11 1.050e-04
3145728 51.9666 45  4.608e-05 1.041e-06  3.2714  4 1.472e-05 20.8572  23 23.1839 11 1.314e-05
}\errorsdB
\pgfplotstableread{
 cells    Gtime  Gits    hdgL2 hdgL2post  qTime  qIt     FEL2  trTime trIt  dgTim dgIt dgL2
    12  0.0020   10  6.842e-01 1.065e+00 0.0014  3  5.621e-01  0.0002  1  0.0023  3 4.684e-01
    96  0.0084   14  1.703e-01 7.177e-02 0.0013  5  1.116e-01  0.0020 10  0.0101 12 9.894e-02
   768  0.0567   17  2.364e-02 4.208e-03 0.0036  5  8.657e-03  0.0098 16  0.0184 11 7.977e-03
  6144  0.6151   22  1.838e-03 1.417e-04 0.0132  5  6.786e-04  0.1558 18  0.0789 11 6.475e-04
 49152  6.8369   31  1.423e-04 5.306e-06 0.1749  6  4.312e-05  1.4048 22  0.6825 10 4.227e-05
393216 82.4786   46  1.012e-05 2.811e-07 1.6487  6  2.707e-06 11.6826 22  6.0616 10 2.685e-06
}\errorsdC
\pgfplotstableread{
 cells    Gtime  Gits   hdgL2 hdgL2post  qTime  qIt    FEL2  trTime trIt  dgTim dgIt dgL2
    12   0.0023  11 3.301e-01 2.892e-01 0.0007  3  3.182e-01  0.0003  1  0.0033  3 2.650e-01
    96   0.0180  15 7.747e-02 2.530e-02 0.0024  6  2.989e-02  0.0053 12  0.0196 13 2.719e-02
   768   0.1497  18 4.413e-03 9.283e-04 0.0064  6  1.741e-03  0.0396 18  0.0359 12 1.670e-03
  6144   1.2825  20 2.233e-04 1.437e-05 0.0405  7  5.232e-05  0.5228 21  0.1868 12 4.994e-05
 49152  12.2677  23 8.568e-06 4.795e-07 0.4689  7  1.703e-06  4.8130 24  1.6184 11 1.607e-06
393216 128.4195  30 3.193e-07 5.890e-09 3.4362  6  5.401e-08 39.5544 24 14.0031 12 5.879e-08
}\errorsdD
\pgfplotstableread{
 cells    Gtime  Gits   hdgL2 hdgL2post  qTime  qIt    FEL2  trTime trIt  dgTim dgIt dgL2
   12  0.0040   12  2.099e-01 2.738e-01 0.0012  3  1.490e-01  0.0006  1  0.0068  3  1.343e-01
   96  0.0398   16  2.975e-02 7.538e-03 0.0043  7  1.144e-02  0.0128 13  0.0447 16  1.064e-02
  768  0.2865   18  1.057e-03 1.342e-04 0.0114  7  2.059e-04  0.1355 20  0.0786 14  1.944e-04
 6144  2.8461   22  2.058e-05 2.419e-06 0.1071  8  4.379e-06  1.3050 22  0.4126 14  4.220e-06
49152 29.6873   27  4.299e-07 1.210e-08 0.9047  7  6.990e-08 11.3176 23  3.3943 13  8.093e-08
}\errorsdE
  \begin{tikzpicture}
    \begin{loglogaxis}[
      title style={at={(1,0.944)},anchor=north east,draw=black,fill=white,font=\footnotesize},
      title={$k=1$},
      width=0.52\textwidth,
      height=0.42\textwidth,
      ylabel={Relative L2 error},
      x label style={at={(0.5,0.05)}},
      tick label style={font=\scriptsize},
      label style={font=\scriptsize},
      legend style={font=\scriptsize},
      xmin=5e-4, xmax=9e1,
      ymin=3e-5, ymax=1e-0,
      grid
      ]
      \addplot table[x={qTime}, y={FEL2}] {\errorsdA};
      \pgfplotsset{cycle list shift=3}
      \addplot table[x={trTime}, y={FEL2}] {\errorsdA};
      \pgfplotsset{cycle list shift=-1}
      \addplot table[x={dgTim}, y={dgL2}] {\errorsdA};
      \pgfplotsset{cycle list shift=0}
      \addplot table[x={Gtime}, y={hdgL2}] {\errorsdA};
    \end{loglogaxis}
  \end{tikzpicture}
  \
  \begin{tikzpicture}
    \begin{loglogaxis}[
      title style={at={(1,0.944)},anchor=north east,draw=black,fill=white,font=\footnotesize},
      title={$k=2$},
      width=0.52\textwidth,
      height=0.42\textwidth,
      x label style={at={(0.5,0.05)}},
      tick label style={font=\scriptsize},
      label style={font=\scriptsize},
      legend style={font=\scriptsize},
      legend to name = legend3dS,
      legend columns = 3,
      xmin=5e-4, xmax=9e1,
      ymin=3e-5, ymax=1e0,
      grid
      ]
      \addplot table[x={qTime}, y={FEL2}] {\errorsdB};
      \addlegendentry{CG mat-free};
      \pgfplotsset{cycle list shift=3}
      \addplot table[x={trTime}, y={FEL2}] {\errorsdB};
      \addlegendentry{CG stat cond};
      \pgfplotsset{cycle list shift=-1}
      \addplot table[x={dgTim}, y={dgL2}] {\errorsdB};
      \addlegendentry{DG-SIP mat-free};
      \pgfplotsset{cycle list shift=0}
      \addplot table[x={Gtime}, y={hdgL2}] {\errorsdB};
      \addlegendentry{HDG trace matrix};
      \addplot[mark=triangle*,dashed,black,every mark/.append style={solid}] table[x={Gtime}, y={hdgL2post}] {\errorsdA};
      \addlegendentry{HDG trace matrix post $k-1$};
    \end{loglogaxis}
  \end{tikzpicture}
  \\
  \begin{tikzpicture}
    \begin{loglogaxis}[
      title style={at={(1,0.944)},anchor=north east,draw=black,fill=white,font=\footnotesize},
      title={$k=3$},
      width=0.52\textwidth,
      height=0.42\textwidth,
      xlabel={Compute time [s]},
      ylabel={Relative L2 error},
      x label style={at={(0.5,0.05)}},
      tick label style={font=\scriptsize},
      label style={font=\scriptsize},
      legend style={font=\scriptsize},
      xmin=5e-4, xmax=9e1,
      ymin=1e-7, ymax=1e-1,
      grid
      ]
      \addplot table[x={qTime}, y={FEL2}] {\errorsdC};
      \pgfplotsset{cycle list shift=3}
      \addplot table[x={trTime}, y={FEL2}] {\errorsdC};
      \pgfplotsset{cycle list shift=-1}
      \addplot table[x={dgTim}, y={dgL2}] {\errorsdC};
      \pgfplotsset{cycle list shift=0}
      \addplot table[x={Gtime}, y={hdgL2}] {\errorsdC};
      \addplot[mark=triangle*,dashed,black,every mark/.append style={solid}] table[x={Gtime}, y={hdgL2post}] {\errorsdB};
    \end{loglogaxis}
  \end{tikzpicture}
  \
  \begin{tikzpicture}
    \begin{loglogaxis}[
      title style={at={(1,0.944)},anchor=north east,draw=black,fill=white,font=\footnotesize},
      title={$k=4$},
      width=0.52\textwidth,
      height=0.42\textwidth,
      xlabel={Compute time [s]},
      x label style={at={(0.5,0.05)}},
      tick label style={font=\scriptsize},
      label style={font=\scriptsize},
      legend style={font=\scriptsize},
      xmin=5e-4, xmax=9e1,
      ymin=1e-7, ymax=1e-1,
      grid
      ]
      \addplot table[x={qTime}, y={FEL2}] {\errorsdD};
      \pgfplotsset{cycle list shift=3}
      \addplot table[x={trTime}, y={FEL2}] {\errorsdD};
      \pgfplotsset{cycle list shift=-1}
      \addplot table[x={dgTim}, y={dgL2}] {\errorsdD};
      \pgfplotsset{cycle list shift=0}
      \addplot table[x={Gtime}, y={hdgL2}] {\errorsdD};
      \addplot[mark=triangle*,dashed,black,every mark/.append style={solid}] table[x={Gtime}, y={hdgL2post}] {\errorsdC};
    \end{loglogaxis}
  \end{tikzpicture}
\\
\ref{legend3dS}
\caption{Accuracy over solver time for 3D Laplacian on the sphere with
  high-order curved boundaries and variable coefficients according to
  Eq.~\eqref{eq:diffusivity_shell}.}
\label{fig:accuracy_3d_shell}
\end{figure}
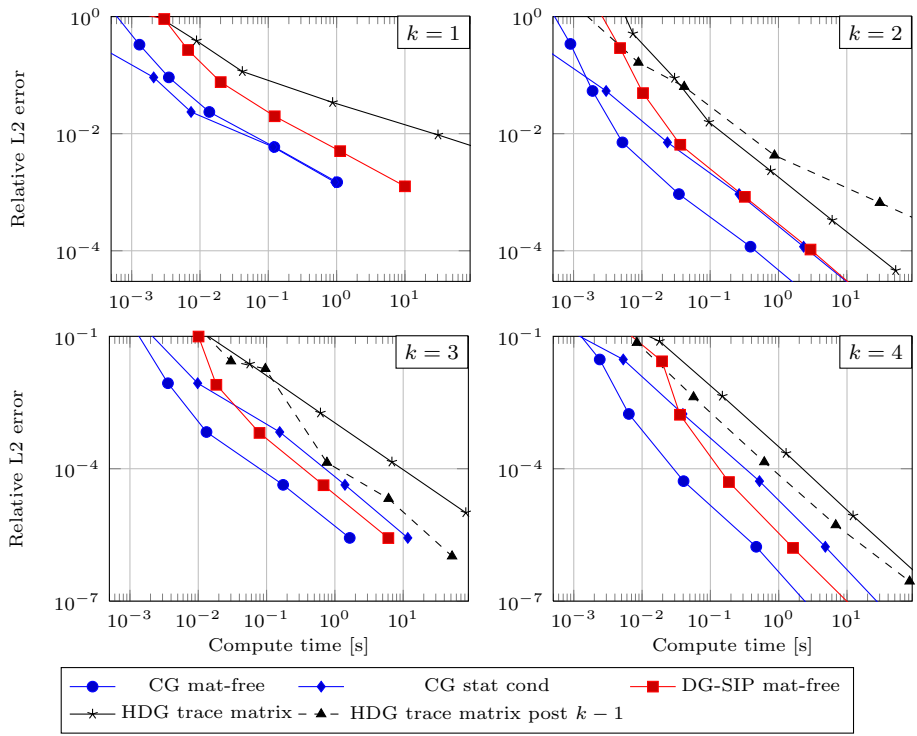

Figs.~\ref{fig:accuracy_3d_cube} and \ref{fig:accuracy_3d_shell} list the
accuracy over the computational time for the constant-coefficient case with
Cartesian mesh and the variable-coefficient case with curved mesh and
high-order mappings. Results appearing in the lower left corner of these plots
combine high accuracy with low computational time. The numbers are from
experiments on a full node with 28 cores of Intel Xeon E5-2690 v4 Broadwell
CPUs. Continuous finite elements show the best efficiency over the whole range
of polynomial degrees. The next best method for $2\leq k \leq 4$, the DG-SIP
method, is two to six times less efficient. The best
matrix-based schemes is the statically condensed finite element method,
slightly ahead of the HDG method with post-processing. Both results are
approximately three to five times slower than DG-SIP and around 20 times
slower than the matrix-free CG implementation for $k\geq 3$. All results are
along the optimal convergence rate curves at $\mathcal O(h^{k+1})$ for the
primal solution $u_h$ and $\mathcal O(h^{k+2})$ for the post-processed
solution in HDG, including the variable-coefficient curved mesh cases.

Fig.~\ref{fig:accuracy_3d_cube} also includes HDG results with polynomials of
complete degree up to $k$, $\mathcal P_k$, spanned by orthogonal Legendre
polynomials, rather than tensor-product space $\mathcal Q_k$. This space skips
the higher order mixed terms such as $xy, xz, yz, xyz$ in $\mathcal Q_1$
elements and thus tightly selects polynomials exactly up to degree $k$. The
$\mathcal P_k$ basis reduces both the number of unknowns by
going from $(k+1)^2$ polynomials per face to $(k+1)(k+2)/2$ and also the nonzero
entries per row, reducing the cost per element by up to a factor of four. The
computational results in Fig.~\ref{fig:accuracy_3d_cube} show that the
decrease in solver times for the tight polynomial space comes with a decrease
in solution accuracy: Even though the solution still converges optimally at
rate $k+1$, the error constants are higher and more elements are needed for
reaching the same accuracy. Thus, no savings can be achieved this way. This
observation is in line with results e.g.~in \cite{ksmw15,Yakovlev16} comparing
tetrahedral to hexahedral element shapes, where a similar or slightly better
efficiency per degree of freedom of hexahedral elements was demonstrated
in the context of matrix-based HDG.

Table~\ref{tab:solver_times_3d} lists the number of iterations to reduce the
linear residual by $10^{-9}$ with the preconditioned conjugate gradient method
as well as solver throughput. Up to 16 million DoFs per second
can be processed with the matrix-free continuous finite element
implementation on 28 cores on Cartesian meshes, and 7.75 million unknowns
per second on a curved mesh. Throughout $2\leq k \leq 6$, the
throughput is above 13 million unknowns per second on the Cartesian
mesh and 6 million unknowns per second on the curved mesh. For comparison, we measured
6.86 million DoFs with HPGMG and polynomial degree 2 on the same
system, again slower than our implementation. This is despite a coarser
iteration tolerance with fewer iterations of HPGMG that relies on a full multigrid
cycle rather than a V-cycle that further reduces the number of operations on
the finest level \cite{hpgmg16}. As reported in \cite{Sundar15}, there is a
slight increase in iteration counts as the polynomial degree $k$ increases,
but high-order methods appear highly attractive nonetheless.

Table~\ref{tab:solver_times_3d} confirms that DG-SIP provides the second
highest solver throughput. Note that the cost per degree of freedom is almost
independent of the polynomial degree for the matrix-free multigrid solvers,
confirming the results of matrix-vector products in
Sec.~\ref{sec:implementation}. The factor between the fastest realization (at
degree between two and four) and the slowest one is only approximately
two. For the matrix-based HDG and statically condensed CG methods, we notice a
distinct decrease in throughput as the polynomial degree increases. This is a
direct consequence of the fact that the matrix rows are more densely populated
for higher degrees, with the cost being directly proportional to this
number. Even though the number of unknowns goes down with static
condensation as compared to matrix-free evaluation,
Fig.~\ref{fig:matvec} shows that this reduction is not nearly enough for
competitive performance. In other words, the gap
between the matrix-free implementations and the matrix-based schemes widens as
the degree increases.

\begin{table}
  \caption{Number of iterations for matrix sizes close to ten million degrees of freedom as well as absolute performance in terms of degrees of freedom solved per second on 28 Broadwell cores. Conjugate gradient tolerances: $10^{-9}$.}
\label{tab:solver_times_3d}
\footnotesize
\centering
\begin{tabular}{lcccccccc}
\hline
 & \multicolumn{2}{c}{CG mat-free} & \multicolumn{2}{c}{CG stat cond} & \multicolumn{2}{c}{DG-SIP mat-free} & \multicolumn{2}{c}{HDG trace matrix} \\
$k$ & its & DoFs/s & its & DoFs/s  & its & DoFs/s & its & DoFs/s \\
\hline
\multicolumn{9}{c}{3D Cartesian mesh, constant coefficients}\\
\hline
1 & 4 & $7.89\cdot 10^6$ & 11 & $3.94\cdot 10^6$ & 14 & $3.07\cdot 10^6$ & 16 & $2.08\cdot 10^6$\\
2 & 5 & $1.35\cdot 10^7$ & 21 & $1.17\cdot 10^6$ & 12 & $5.17\cdot 10^6$ & 15 & $1.14\cdot 10^6$\\
3 & 5 & $1.45\cdot 10^7$ & 22 & $6.17\cdot 10^5$ & 10 & $6.42\cdot 10^6$ & 17 & $6.15\cdot 10^5$\\
4 & 5 & $1.58\cdot 10^7$ & 25 & $3.46\cdot 10^5$ & 10 & $6.33\cdot 10^6$ & 17 & $4.11\cdot 10^5$\\
5 & 5 & $1.52\cdot 10^7$ & 26 & $2.29\cdot 10^5$ & 11 & $5.03\cdot 10^6$ & 20 & $2.48\cdot 10^5$\\
6 & 5 & $1.43\cdot 10^7$ & 27 & $1.98\cdot 10^6$ & 11 & $5.06\cdot 10^6$ & 20 & $1.87\cdot 10^5$\\
7 & 5 & $1.25\cdot 10^7$ & 26 & $1.30\cdot 10^5$ & 13 & $4.04\cdot 10^6$ & 21 & $1.40\cdot 10^5$\\
8 & 5 & $1.17\cdot 10^7$ & 27 & $9.90\cdot 10^4$ & 12 & $3.64\cdot 10^6$ & 21 & $1.13\cdot 10^5$\\
\hline
\multicolumn{9}{c}{3D curved mesh, variable coefficients}\\
\hline
1 & 5 & $3.13\cdot 10^6$ & 14 & $3.38\cdot 10^6$ & 13 & $2.51\cdot 10^6$ & 217 & $1.56\cdot 10^5$\\
2 & 5 & $7.75\cdot 10^6$ & 23 & $1.06\cdot 10^6$ & 11 & $3.36\cdot 10^6$ & 32 & $5.43\cdot 10^5$\\
3 & 6 & $6.51\cdot 10^6$ & 22 & $6.40\cdot 10^5$ & 10 & $4.15\cdot 10^6$ & 46 & $2.31\cdot 10^5$\\
4 & 5 & $7.38\cdot 10^6$ & 24 & $3.68\cdot 10^5$ & 12 & $3.51\cdot 10^6$ & 30 & $2.32\cdot 10^5$\\
5 & 7 & $6.88\cdot 10^6$ & 23 & $2.65\cdot 10^5$ & 13 & $3.12\cdot 10^6$ & 27 & $1.82\cdot 10^5$\\
6 & 8 & $6.08\cdot 10^6$ & 24 & $1.85\cdot 10^5$ & 15 & $2.65\cdot 10^6$ & 24 & $1.57\cdot 10^5$\\
7 & 10 & $4.90\cdot 10^6$ & 25 & $1.35\cdot 10^5$ & 17 & $2.35\cdot 10^6$ & 27 & $1.09\cdot 10^5$\\
8 & 11 & $4.47\cdot 10^6$ & 25 & $1.12\cdot 10^5$ & 21 & $1.95\cdot 10^6$ & 22 & $1.14\cdot 10^5$\\
\hline
\end{tabular}
\end{table}

Note that the smaller iteration counts for the geometric multigrid approaches
are due to a more expensive smoother that involves five matrix-vector
products rather than only one forward and backward substitution in the ILU of
HDG. The HDG solver with $p$-AMG is highly competitive in terms of the total
number of operator evaluations to reach the prescribed tolerance of
$10^{-9}$. For example, the HDG solver involves 45 operator evaluations and 30
ILU applications on the finest level in the Cartesian case with $k=2$ of
Table~\ref{tab:solver_times_3d} that lists 15 iterations, as compared to 65
matrix-vector products for the CG solver with GMG preconditioning at 5
iterations and 156 matrix-vector products for DG-SIP at 12 iterations. It is
rather the different performance of matrix-vector products that favors the
matrix-free schemes.

\begin{figure}
\pgfplotstableread{
degree  sizefe   feqmf   feqsc  sizedg    dgsipmf  sizehdg   hdg    sizehdgpost hdgpost
1       37761025 3.6890  6.3680 150994944 29.16    205520896 207.5  1640960     0.8499
2       591361   0.0185  0.1169 1327104   0.0979   1575936   0.6617 99072       0.0163
3       83521    0.0048  0.0092 147456    0.0151   132096    0.0321 33280       0.0060
4       37249    0.0036  0.0051 57600     0.0095   31920     0.0079 6000        0.0033
5       14641    0.0032  0.0035 20736     0.0114   9744      0.0049 3264        0.0030
6       5329     0.0033  0.0008 12544     0.0057   5880      0.0039 2184        0.0030
}\mytimestwo
\pgfplotstableread{
degree  sizefe     feqmf     feqsc     sizedg     dgsipmf   sizehdg    hdg      sizehdgpost hdgpost
1       3.7325e+11 4.52e+04  5.816e+04 5.9463e+12 4.67e+05  4.44e+12   2.84e+06 3.1654e+09  916.53
2       454756609  32.2      289.15    1.5288e+09 296.92    3.6239e+09 3585.2   38271744    35.2064
3       38272753   2.2760    40.085    89915392   14.274    100663296  166.36   5419008     8.7499
4       7189057    0.4315    11.6710   13824000   2.1823    13406400   32.6582  1705200     4.0742
5       2803221    0.1582    5.8964    2985984    0.5079    3649536    13.9873  907200      3.3229
6       912673     0.0441    1.9536    1404928    0.2723    1234800    5.9680   275184      1.2823
7       614125     0.0469    1.5820    884736     0.3409    835584     5.6285   359424      2.2516
8       912673     0.0684    2.7767    1259712    0.5743    454896     3.7381   139968      1.1430
}\mytimes
\centering
  \begin{tikzpicture}
    \begin{semilogyaxis}[
      title style={at={(1,0.943)},anchor=north east,draw=black,fill=white,font=\footnotesize},
      title={2D Cartesian mesh},
      width=0.48\textwidth,
      height=0.42\textwidth,
      xlabel={Polynomial degree $k$},
      ylabel={Compute time [s]},
      x label style={at={(0.5,0.05)}},
      tick label style={font=\scriptsize},
      label style={font=\scriptsize},
      legend style={font=\scriptsize},
      legend to name = legendTimeAcc,
      legend columns = 3,
      ymin=1e-3, ymax=5e2,
      grid
      ]
      \addplot table[x={degree}, y={feqmf}] {\mytimestwo};
      \addlegendentry{CG mat-free};
      \pgfplotsset{cycle list shift=3}
      \addplot table[x={degree}, y={feqsc}] {\mytimestwo};
      \addlegendentry{CG stat cond};
      \pgfplotsset{cycle list shift=-1}
      \addplot table[x={degree}, y={dgsipmf}] {\mytimestwo};
      \addlegendentry{DG-SIP mat-free};
      \pgfplotsset{cycle list shift=0}
      \addplot table[x={degree}, y={hdg}] {\mytimestwo};
      \addlegendentry{HDG trace matrix};
      \addplot[mark=triangle*,dashed,black,every mark/.append style={solid}] table[x expr=\thisrowno{0}, y index=9] {\mytimestwo};
      \addlegendentry{HDG trace matrix post};
    \end{semilogyaxis}
  \end{tikzpicture}
  \begin{tikzpicture}
    \begin{semilogyaxis}[
      title style={at={(1,0.943)},anchor=north east,draw=black,fill=white,font=\footnotesize},
      title={3D Cartesian mesh},
      width=0.54\textwidth,
      height=0.42\textwidth,
      xlabel={Polynomial degree $k$},
      x label style={at={(0.5,0.05)}},
      tick label style={font=\scriptsize},
      ytick={1e-2,1e-1,1,1e1,1e2,1e3,1e4},
      label style={font=\scriptsize},
      legend style={font=\scriptsize},
      ymin=2e-2, ymax=9e4,
      grid
      ]
      \addplot table[x={degree}, y={feqmf}] {\mytimes};
      \pgfplotsset{cycle list shift=3}
      \addplot table[x={degree}, y={feqsc}] {\mytimes};
      \pgfplotsset{cycle list shift=-1}
      \addplot table[x={degree}, y={dgsipmf}] {\mytimes};
      \pgfplotsset{cycle list shift=0}
      \addplot table[x={degree}, y={hdg}] {\mytimes};
      \addplot[mark=triangle*,dashed,black,every mark/.append style={solid}] table[x expr=\thisrowno{0}, y index=9] {\mytimes};
    \end{semilogyaxis}
  \end{tikzpicture}
  \\
  \ref{legendTimeAcc}
  \caption{Time to reach a discretization accuracy of $10^{-6}$ as a function
    of the polynomial degree.}
\label{fig:time_accuracy_3d}
\end{figure}
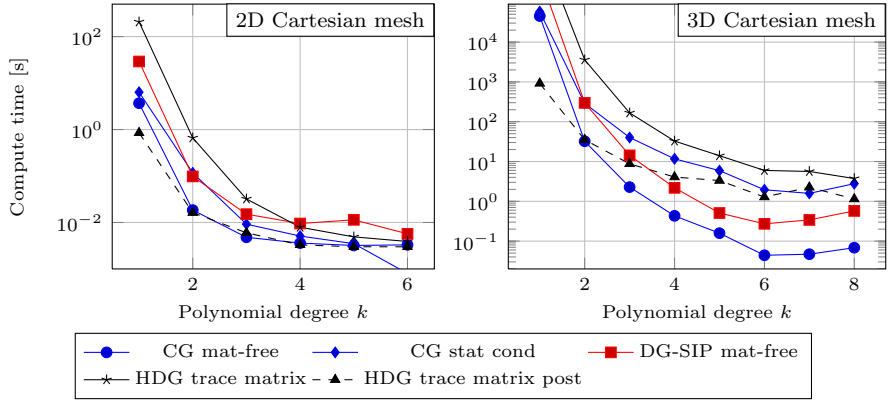

Fig.~\ref{fig:time_accuracy_3d} lists the computational time required to reach
a fixed relative accuracy of $10^{-6}$ for various polynomial degrees on a
Cartesian mesh with 28 Broadwell cores, both in two and three space
dimensions. Note that the 3D numbers for $k=1$ need significantly more memory
and computational resources than what is available the (fat-memory) single
node with 512 GB used for the present tests, requiring $3.7\cdot 10^{11}$
degrees of freedom for continuous elements or $4.4\cdot 10^{12}$ degrees of
freedom in the trace system. For reasons of comparison, extrapolations of the
computational time recorded at around $10^8$ degrees of freedom to accuracy
$10^{-6}$ have been used under the justified assumption of optimal iteration
counts and convergence rates. The efficiency dramatically increases as the
polynomial degree is risen, enabling the solution to a tolerance of $10^{-6}$
in less than 10 seconds for degree $k=3$ with the matrix-free CG method and
the post-processed HDG method at $k=3$ in three dimensions, or $0.01$ seconds
in 2D.

When increasing the polynomial degree further, different saturation points
appear in 3D. Matrix-based schemes mainly suffer
from the aforementioned increase of nonzero entries per row, despite the
number of DoFs still going down, a behavior also described in
\cite{Gholami16}. The faster matrix-free schemes with approximately constant timings
per unknown run into latency issues instead, including the coarse grid
solver. Furthermore, the granularity of the mesh sizes that are tested in the
pattern $4^3$, $8^3$, $12^3$, $16^3$, $20^3$, $24^3$, $28^3$, $32^3$, $40^3$,
$48^3$ (and continuing with multiples of $4,5,6,7$ times a power of two)
result in selecting the mesh $12^3$ for all of $k=6,7,8$, whereas the next
coarser size $8^3$ is too coarse unless $k\geq 9$.

\subsection{Three-dimensional example with non-smooth solution}

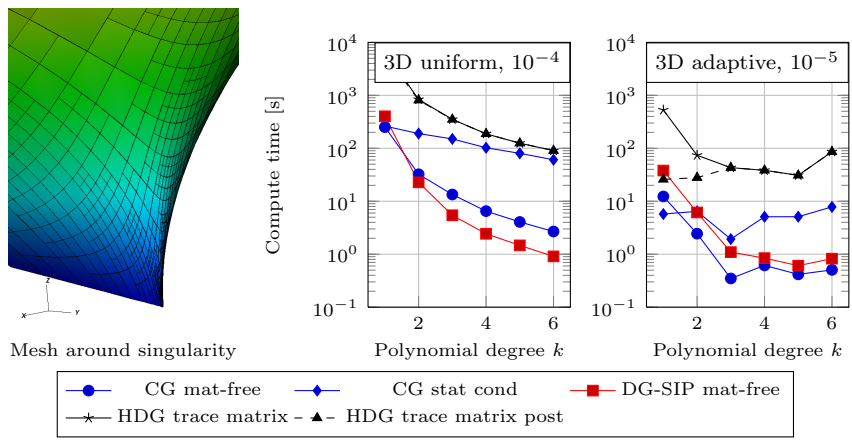
\begin{figure}
\pgfplotstableread{
degree  sizefe     feqmf     feqsc     sizedg     dgsipmf   sizehdg    hdg      sizehdgpost hdgpost
1       1.2198e+09 251.52    262.65    1.4033e+09 402.00    1.1688e+10 6767.6   1.0445e+10  6048.3
2       3.2383e+08 32.357    189.71    1.1466e+08 22.732    1.0378e+09 835.72   1.0052e+09  809.48
3       1.4148e+08 13.412    150.25    3.4618e+07 5.4440    2.4579e+08 352.11   2.4492e+08  350.86
4       7.6933e+07 6.5084    102.40    1.5425e+07 2.4320    8.7822e+07 187.82   8.7822e+07  187.82
5       4.5886e+07 4.0638    80.270    8.1247e+06 1.4699    3.9076e+07 124.82   3.9076e+07  124.82
6       2.9688e+07 2.6841    60.760    4.9193e+06 0.90986   2.0076e+07 90.595   2.0076e+07  90.595
}\timesnonsmoothuni
\pgfplotstableread{
degree  sizefe     feqmf     feqsc     sizedg     dgsipmf   sizehdg    hdg      sizehdgpost hdgpost
1       24454468   12.324    5.7356    88169040   38.024    69881300   528.85   137214320   25.875
2       15741814   2.4426    6.4486    20841138   6.1455    24958377   73.8373  11706489    27.8722
3       3085649    0.3482    1.9304    6094848    1.0953    11664288   43.1149  11664288    43.1149
4       5333507    0.6147    5.0910    4197000    0.8433    6532200    38.3779  6532200     38.3779
5       4069979    0.4162    5.1082    2722464    0.6082    3798000    30.8705  3798000     30.8705
6       4794119    0.5064    7.8297    2949800    0.8243    7164682    85.932   7164682     85.932
}\timesnonsmoothada
\centering
\begin{tikzpicture}
\node at (0,0){\includegraphics[width=0.25\textwidth]{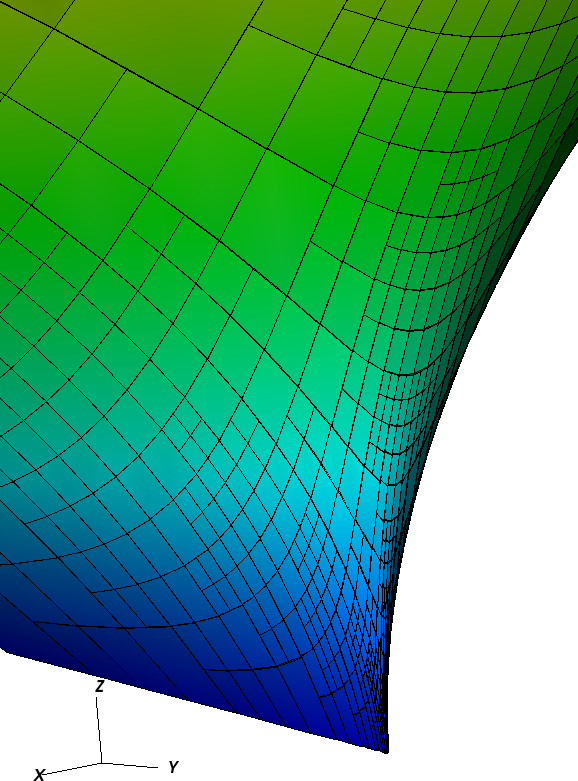}};
\node at (0,-2.5){\scriptsize Mesh around singularity};
\end{tikzpicture}
  \begin{tikzpicture}
    \begin{semilogyaxis}[
      title style={at={(1,0.94)},anchor=north east,draw=black,fill=white,font=\footnotesize},
      title={3D uniform, $10^{-4}$},
      width=0.35\textwidth,
      height=0.42\textwidth,
      xlabel={Polynomial degree $k$},
      ylabel={Compute time [s]},
      x label style={at={(0.5,0.05)}},
      y label style={at={(0.,0.5)}},
      tick label style={font=\scriptsize},
      ytick={1e-2,1e-1,1,1e1,1e2,1e3,1e4},
      label style={font=\scriptsize},
      legend style={font=\scriptsize},
      legend to name = legendTimeAccNon,
      legend columns = 3,
      ymin=1e-1, ymax=1e4,
      grid
      ]
      \addplot table[x={degree}, y={feqmf}] {\timesnonsmoothuni};
      \addlegendentry{CG mat-free};
      \pgfplotsset{cycle list shift=3}
      \addplot table[x={degree}, y={feqsc}] {\timesnonsmoothuni};
      \addlegendentry{CG stat cond};
      \pgfplotsset{cycle list shift=-1}
      \addplot table[x={degree}, y={dgsipmf}] {\timesnonsmoothuni};
      \addlegendentry{DG-SIP mat-free};
      \pgfplotsset{cycle list shift=0}
      \addplot table[x={degree}, y={hdg}] {\timesnonsmoothuni};
      \addlegendentry{HDG trace matrix};
      \addplot[mark=triangle*,dashed,black,every mark/.append style={solid}] table[x expr=\thisrowno{0}, y index=9] {\timesnonsmoothuni};
      \addlegendentry{HDG trace matrix post};
    \end{semilogyaxis}
  \end{tikzpicture}
  \begin{tikzpicture}
    \begin{semilogyaxis}[
      title style={at={(1,0.94)},anchor=north east,draw=black,fill=white,font=\footnotesize},
      title={3D adaptive, $10^{-5}$},
      width=0.35\textwidth,
      height=0.42\textwidth,
      xlabel={Polynomial degree $k$},
      x label style={at={(0.5,0.05)}},
      tick label style={font=\scriptsize},
      ytick={1e-2,1e-1,1,1e1,1e2,1e3,1e4},
      label style={font=\scriptsize},
      legend style={font=\scriptsize},
      ymin=1e-1, ymax=1e4,
      grid
      ]
      \addplot table[x={degree}, y={feqmf}] {\timesnonsmoothada};
      \pgfplotsset{cycle list shift=3}
      \addplot table[x={degree}, y={feqsc}] {\timesnonsmoothada};
      \pgfplotsset{cycle list shift=-1}
      \addplot table[x={degree}, y={dgsipmf}] {\timesnonsmoothada};
      \pgfplotsset{cycle list shift=0}
      \addplot table[x={degree}, y={hdg}] {\timesnonsmoothada};
      \addplot[mark=triangle*,dashed,black,every mark/.append style={solid}] table[x expr=\thisrowno{0}, y index=9] {\timesnonsmoothada};
    \end{semilogyaxis}
  \end{tikzpicture}
  \\
  \ref{legendTimeAccNon}
  \caption{Singular solution. Time to reach a discretization accuracy of
    $10^{-4}$ on the uniform mesh and $10^{-5}$ on the adaptive mesh as a
    function of the polynomial degree on 28 cores.}
\label{fig:time_accuracy_nonsmooth}
\end{figure}

We now consider the Laplacian on a cube $(-1,1)^3$ with a slit along the plane
$\{x=0, y<0, -1<z<1\}$. The solution is given by
\[
u(x,y,z) = r^{\frac 12} + \sin(0.5 \phi)\quad \text{with}\quad  r=\sqrt{x^2+y^2},
\phi = \arctan\left(\frac{y}{x}\right)+\pi,
\]
and constant in $z$-direction. Fig.~\ref{fig:time_accuracy_nonsmooth} shows
the solution around the singularity on a slice at $z=0$ with elevation along
the function value. Due to the singularity, convergence rates in the $L_2$
norm are only linear in the mesh size, irrespective the polynomial
degree. Fig.~\ref{fig:time_accuracy_nonsmooth} lists the time to reach a
discretization accuracy of $10^{-4}$ on uniform meshes and $10^{-5}$ with
adaptive meshes which are created by successively refining the 15\% of the
elements with the largest jump in the gradient over element boundaries as a
simple error estimator \cite{Gago83}. The results confirm the considerably
higher efficiency of the matrix-free schemes both in the uniform and adaptive
mesh case. Furthermore, the methods become more efficient as the polynomial
degree is increased also on the adaptive mesh, as opposed to the matrix-based
methods that level off and high order methods do not pay off. Also note that
DG-SIP is more efficient than continuous elements on the uniform mesh due to
better solution accuracy around the singularity, allowing for coarser meshes.

\subsection{Scalability in the massively parallel context}

\begin{figure}
\pgfplotstableread{
nprocs    dg256k       fem256k    dg2m     fem2m    dg16m    fem16m   dg128m   fe128m   hdg256k hdg2m   hdg16m
32        1.3247       0.325741   10.8376  2.51771  nan      nan      nan      nan      nan     nan     nan
64        0.6587       0.1645     5.4530   1.2823   nan      nan      nan      nan      7.57474 nan     nan
128       0.330609     0.090898   2.76268  0.658832 nan      5.00366  nan      nan      3.7768  nan     nan
256       0.1746       0.059922   1.3686   0.339999 11.2383  2.56216  nan      nan      1.91563 nan     nan
512       0.103789     0.0455449  0.684803 0.176482 5.66666  1.29986  nan      nan      0.98943 8.07101 nan
1024      0.0685959    0.0368049  0.345644 0.099691 2.926    0.67364  nan      6.48155  0.5803  4.0690  nan
2048      0.0524418    0.0348921  0.184925 0.069573 1.48659  0.356066 12.5784  2.59601  0.2697  2.2784  nan
4096      0.045135     0.0367949  0.114661 0.056833 0.758113 0.19293  5.85988  1.3251   0.1657  1.2719  12.387
8192      0.0443981    0.033958   0.075846 0.045350 0.405509 0.110485 3.05784  0.790214 0.2517  0.92450 6.529
16384     0.0411382    0.0379629  0.06585  0.049351 0.251308 0.099904 1.58422  0.424692 0.2429  0.59247 3.6394
32768     0.049526     0.0461671  0.056096 0.051276 0.155227 0.077546 0.860195 0.229114 0.2652  0.57048 1.9543
65536     0.0588942    0.0466189  0.070462 0.058941 0.122254 0.075194 0.480742 0.157353 0.2343  0.59921 1.2854
147456    0.0618479    nan        0.0699389 nan     0.1052   nan      0.289767 nan      nan     nan     nan
}\scalingSNB
\pgfplotstableread{
nprocs    dg256k       fem256k    dg2m     fem2m    dg16m    fem16m   dg128m   fe128m   hdg256k hdg2m   hdg16m
28        1.1358       0.3501     9.0596   2.8813   nan      nan      nan      nan      nan     nan     nan
56        0.567158     0.1610     4.61779  1.4397   nan      nan      nan      nan      nan     nan     nan
112       0.270517     0.0785     2.37399  0.7235   nan      5.7292   nan      nan      4.95074 nan     nan
224       0.1299       0.0441     1.1955   0.3665   9.3419   2.9128   nan      nan      2.62977 nan     nan
448       0.0816212    0.0265     0.611527 0.1853   4.80802  1.4754   nan      nan      1.7138  11.8409 nan
896       0.0468       0.0256     0.3175   0.0920   2.4471   0.7511   nan      5.7863   0.9265  6.12811 nan
1792      0.0411       0.0226     0.1825   0.0571   1.2266   0.4088   9.5697   2.9563   0.5768  3.0865  nan
3584      0.0301       0.0202     0.0986   0.0447   0.6498   0.2211   4.8884   1.5270   0.3943  1.6857  13.387
7168      0.0232       0.0174     0.0720   0.0337   0.3802   0.1395   2.5283   0.7975   0.3717  0.92450 8.026
14336     0.0197902    0.0178611  0.0455   0.025414 0.2167   0.084578 1.3324   0.449233 0.3729  0.59247 4.93898
}\scalingHSW
\centering
  \begin{tikzpicture}
    \begin{loglogaxis}[
      title style={at={(1,0.955)},anchor=north east,draw=black,fill=white,font=\footnotesize},
      title={$128^3$ mesh, $\mathcal Q_3$ elements},
      width=0.53\textwidth,
      height=0.5\textwidth,
      xlabel={Number of cores},
      ylabel={Solver time [s]},
      x label style={at={(0.5,0.05)}},
      y label style={at={(0.05,0.5)}},
      xtick={32,128,512,2048,8192,32768,131072},
      xticklabels={32,128,512,2048,8192,32k,128k},
      tick label style={font=\scriptsize},
      label style={font=\scriptsize},
      legend style={font=\scriptsize},
      legend to name = legendstrong,
      legend columns = 3,
      ymin=2e-2, ymax=30,
      grid
      ]
      \addplot table[x={nprocs}, y={fem2m}] {\scalingSNB};
      \addlegendentry{CG SandyB};
      \addplot[blue,mark=o] table[x={nprocs}, y={fem2m}] {\scalingHSW};
      \addlegendentry{CG Haswell};
      \pgfplotsset{cycle list shift=-1}
      \addplot table[x={nprocs}, y={dg2m}] {\scalingSNB};
      \addlegendentry{DG-SIP SandyB};
      \addplot[red,mark=square] table[x={nprocs}, y={dg2m}] {\scalingHSW};
      \addlegendentry{DG-SIP Haswell};
      \addplot table[x={nprocs}, y={hdg2m}] {\scalingSNB};
      \addlegendentry{HDG SandyB};
      \pgfplotsset{cycle list shift=0}
      \addplot[dashed,black] coordinates {
        (16,10)
        (131072,10/8192)
      };
      \addlegendentry{linear scaling};
    \end{loglogaxis}
  \end{tikzpicture}
  \
  \begin{tikzpicture}
    \begin{loglogaxis}[
      title style={at={(1,0.955)},anchor=north east,draw=black,fill=white,font=\footnotesize},
      title={$256^3$ mesh, $\mathcal Q_3$ elements},
      width=0.53\textwidth,
      height=0.5\textwidth,
      xlabel={Number of cores},
      x label style={at={(0.5,0.05)}},
      xtick={32,128,512,2048,8192,32768,131072},
      xticklabels={32,128,512,2048,8192,32k,128k},
      tick label style={font=\scriptsize},
      label style={font=\scriptsize},
      legend style={font=\scriptsize},
      ymin=2e-2, ymax=30,
      grid
      ]
      \addplot table[x={nprocs}, y={fem16m}] {\scalingSNB};
      \addplot[blue,mark=o] table[x={nprocs}, y={fem16m}] {\scalingHSW};
      \pgfplotsset{cycle list shift=-1}
      \addplot table[x={nprocs}, y={dg16m}] {\scalingSNB};
      \addplot[red,mark=square] table[x={nprocs}, y={dg16m}] {\scalingHSW};
      \addplot table[x={nprocs}, y={hdg16m}] {\scalingSNB};
      \addplot[dashed,black] coordinates {
        (128,10)
        (131072,10/1024)
      };
    \end{loglogaxis}
  \end{tikzpicture}
  \\
  \ref{legendstrong}
  \caption{Strong scaling experiment on SuperMUC Phase 1 (SandyB) and SuperMUC
    Phase 2 (Haswell) on up to 147\,456 cores.}
\label{fig:strong_scaling}
\end{figure}
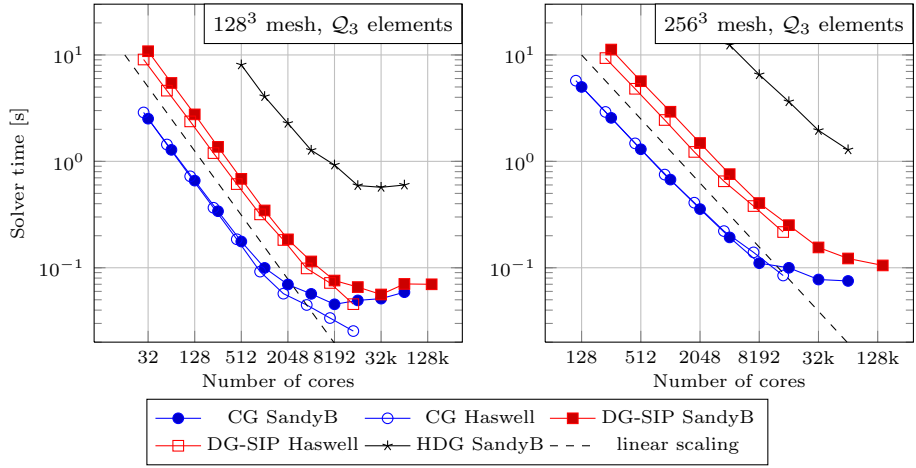

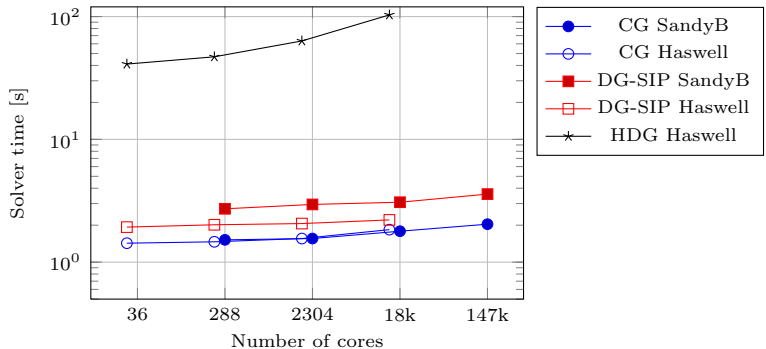
\begin{figure}
\centering
  \begin{tikzpicture}
    \begin{loglogaxis}[
      width=0.6\textwidth,
      height=0.45\textwidth,
      xlabel={Number of cores},
      ylabel={Solver time [s]},
      x label style={at={(0.5,0.05)}},
      y label style={at={(0.05,0.5)}},
      xtick={36,288,2304,18432,147456},
      xticklabels={36,288,2304,18k,147k},
      tick label style={font=\scriptsize},
      label style={font=\scriptsize},
      legend style={font=\scriptsize},
      legend pos = outer north east,
      ymin=5e-1, ymax=120,
      grid
      ]
      \addplot table[x index=0, y expr={\thisrowno{0}*1e6*\thisrowno{2}/\thisrowno{1}/27}] {
nprocs  nele    time     nele   time
288     7.08e6  1.00813  16.8e6 2.33326
2304    56.6e6  1.03276  134e6  2.35964
18432   453e6   1.18549  1.07e9 2.56162
147456  3.62e9  1.35056  8.6e9  2.90856
      };
      \addlegendentry{CG SandyB};
      \addplot[blue,mark=o] table[x index=0, y expr={\thisrowno{0}*1e6*\thisrowno{6}/\thisrowno{5}/27}] {
  procs     nele     time      nele   time      nele     time
  28        111e3    0.130938  262e3  0.350081  885e3    1.21754
  224       885e3    0.143288  2.10e6 0.366498  7.078e6  1.24945
  1792      7.078e6  0.204185  16.8e6 0.408801  56.6e6   1.32782
  14336     56.62e6  0.2234    135e6  0.4767    453e6    1.56934
      };
      \addlegendentry{CG Haswell};
      \pgfplotsset{cycle list shift=-1}
      \addplot table[x index=0, y expr={\thisrowno{0}*1e6*\thisrowno{4}/\thisrowno{3}/64}] {
nprocs  nele    time     nele    time     nele   time
288     2.1e6   1.23979  7.08e6  4.2793   16.8e6 11.2383
2304    16.8e6  1.37338  56.6e6  4.63946  134e6  11.2555
18432   134e6   1.53388  453e6   4.83723  1.07e9 11.9156
147456  1.07e9  1.6509   3.62e9  5.6287   8.6e9  13.4813
      };
      \addlegendentry{DG-SIP SandyB};
      \addplot[red,mark=square] table[x index=0, y expr={\thisrowno{0}*1e6*\thisrowno{6}/\thisrowno{5}/64}] {
  procs     nele     time      nele   time      nele     time
  28        111e3    0.466557  262e3  1.13582  885e3    3.90259
  224       885e3    0.514658  2.10e6 1.19551  7.078e6  4.07576
  1792      7.078e6  0.581852  16.8e6 1.22658  56.6e6   4.1678
  14336     56.62e6  0.6202    135e6  1.3324   453e6    4.4648
      };
      \addlegendentry{DG-SIP Haswell};
      \addplot table[x index=0, y expr={\thisrowno{0}*1e6*\thisrowno{2}/\thisrowno{1}/48}] {
nprocs    ncells    timehdg
28        64000     4.5033
224       512000    5.1608
1792      4096000   6.9472
14336     32768000  11.323
      };
      \addlegendentry{HDG Haswell};
    \end{loglogaxis}
  \end{tikzpicture}
  \caption{Weak scaling experiment on SuperMUC Phase 1 (SandyB) and SuperMUC
    Phase 2 (Haswell) using a grain size of 1 million unknowns per core and
    $\mathcal Q_3$ elements in 3D.}
\label{fig:weak_scaling}
\end{figure}

The parallel behavior of the solvers is shown in Fig.~\ref{fig:strong_scaling}
by a strong scaling experiment and in Fig.~\ref{fig:weak_scaling} by a weak
scaling experiment. All codes have been parallelized with pure MPI according
to the techniques described in \cite{bbhk11,p4est,kk12}. Two large-scale
parallel systems have been used, SuperMUC Phase 1 consisting of up to 9216
nodes with $2\times 8$ cores (Intel Xeon E5-2680 Sandy Bridge, 2.7 GHz), and
SuperMUC Phase 2 consisting of up to 512 nodes with $2\times 14$ cores (Intel
Xeon E5-2697 v3 Haswell, 2.6 GHz). In Fig.~\ref{fig:strong_scaling}, we
observe ideal strong scaling of the geometric multigrid solvers until a lower
threshold of approximately 0.05--0.1 seconds where communication latency
becomes dominant. We note that the scaling in the GMG solvers for CG and
DG-SIP saturates below 30\,000 degrees of freedom per core in both panels of
Fig.~\ref{fig:strong_scaling}. Also note the wide range of problem sizes with
almost two orders of magnitude going from saturated scaling to the size that
still fits into approximately 2 GB RAM memory per core, much more than for the
matrix-based realization. On the other hand, the HDG solver is already saturated at
around 0.5 seconds. This breakdown is due to the non-ideal behavior of the
ML-AMG part also reported in \cite{Sundar12}. The scaling of the HDG solver is
relatively good until matrix sizes go below 5000 rows per core. Note that on a
smaller $64^3$ mesh, nearly linear scaling down to approximately $0.2\text{s}$
has been obtained.

The weak scaling plots in Fig.~\ref{fig:weak_scaling} display the time to
solve a linear system with one million degrees of freedom per core for the
three chosen discretizations as the number of processors and the problem size
increase at the same rate.  Between the smallest and largest configuration in
the weak scaling tests, parallel efficiencies of 75\%, 78\%, 75\%, and 87\%
have been measured for the CG method on SuperMUC Phase 1 and Phase 2 as well
as DG-SIP on SuperMUC Phase 1 and Phase 2, respectively. For the HDG linear
solvers, the weak scaling is somewhat worse, reaching 40\% when going from 28
to 14\,336 Haswell cores of SuperMUC Phase 2. Approximately half of the
decrease in efficiency is due to the increase in solver iterations from 18 to
29, and the other half is due to inefficiencies in the AMG hierarchy.

The largest computation on 147\,456 cores achieved an arithmetic throughput of
about 1.4 PFLOP/s for DG-SIP, out of a theoretical peak of 3.2 PFLOP/s on
SuperMUC Phase 1. To put the obtained results into perspective, we compare the
CG solution time of $1.52\text{s}$ on 288 cores and $2.04\text{s}$ on 147\,456
cores to the numbers from \cite{Gholami16} which were obtained on the Stampede
system with the same Intel Xeon E5-2680 Sandy Bridge processors, outperforming
both the $2.5\text{s}$ for HPGMG at $\mathcal Q_2$ elements and more than
$20\text{s}$ for GMG implementation of the authors from \cite{Gholami16} on
fourth degree elements. In addition, we could solve a DG-SIP system with one
million degrees of freedom per core in $1.92\text{s}$ on SuperMUC Phase
2. Note that the computational time per core on the Sandy Bridge and Haswell
systems, respectively, is similar for the continuous elements (limited mainly
by indirect addressing into vectors), whereas there is a considerably
advantage of the Haswell system with fused multiply-add instruction and
faster L1 cache access for DG-SIP.

\section{Conclusions}\label{sec:conclusions}

In this study, a performance comparison between continuous and discontinuous
Galerkin methods has been presented. The work has concentrated on
state-of-the-art multigrid solvers for the Laplacian as the prototype elliptic
equation. As opposed to previous studies that focused on direct solvers, our
experiments show that the primal formulation in terms of
continuous finite elements or discontinuous Galerkin symmetric interior
penalty methods allows for up to an order of magnitude more efficient solution than the HDG
method in 3D, also when including superconvergence of HDG. In two space
dimensions, the performance gap is approximately a factor of two to five times
for polynomial degrees $2\leq k\leq 5$. When comparing the continuous finite
element implementation against the interior penalty discontinuous Galerkin
method, a performance advantage of a factor of two to three for the continuous
case has been recorded.

Our results are due to the beneficial properties of modern sum factorization implementations
on quadrilateral and hexahedral meshes. Depending on the structure of the
equations and the behavior of the solution, either continuous finite elements
or symmetric interior penalty discontinuous Galerkin methods are the preferred choice.
We found that the time to solution per degree of freedom is almost
constant for polynomial degrees between two and eight at a similar rate as for the HPGMG benchmark, and better than
matrix-based multigrid schemes on linear finite elements. Thus, the polynomial
degree can be chosen as high as the meshing of the geometry allows for without
compromising throughput. The promising results in this study motivate to
pursue developments of matrix-free solvers in non-elliptic contexts, where the
Jacobi-related techniques used here are not sufficient and matrix-based
methods use Gauss--Seidel or ILU smoothers.

Our conclusions go against the results of previous efficiency studies and are
mainly explained by the different performance of matrix-vector products. In
particular, counting degrees of freedom or arithmetic operations is not enough
to judge application performance. The higher performance is due to a
reduced memory transfer as compared to memory-limited sparse matrix kernels. A
roofline performance model has been developed which shows that our results are
close to the performance limits of the underlying hardware, making the
conclusions general without bias towards any of the methods.
Moreover, the performance model allows for predicting performance on
other HPC systems with different machine balances. Note that reaching the high
performance numbers recorded for the matrix-free solvers needs careful
implementation that is only realistic for large finite element libraries that
can distribute the development burden over many applications. However, we
expect the developments from \cite{hpgmg16,brown10,kk12} to become mainstream
library components of finite element codes in the future, similarly to
high-performance dense linear algebra implementations of BLAS and
LAPACK. In order to exploit the fast matrix-free approaches also in the
context of HDG, alternative evaluation schemes have been presented in this
work. The HDG mixed system involving the primal variable and the flux is
clearly faster than the sparse matrix-vector product with the trace matrix,
despite considerably fewer degrees of freedom in the latter. Given appropriate
solvers, these approaches promise best HDG performance.

\appendix

\section*{Acknowledgments}
The authors would like to thank Katharina Kormann and Niklas Fehn for
discussions about the manuscript and Timo Heister and Guido Kanschat on
the multigrid implementation in \texttt{deal.II}.

\end{document}